\documentclass[12pt,a4]{article}
\usepackage{amssymb,amsmath,graphicx,enumerate}
\usepackage{longtable,mathrsfs}
\usepackage{hyperref}
\graphicspath {P:\Dole}

\DeclareMathAlphabet{\mathpzc}{OT1}{pzc}{m}{it}

\numberwithin{equation}{section}

\newtheorem{theorem}{Theorem}[section]
\newtheorem{lemma}{Lemma}[section]

\newtheorem{definition}{Definition}[section]
\newtheorem{example}{Example}[section]

\title{Canonical sine and cosine Transforms For Integrable Boehmians}

\author{Pravinkumar V. Dole and S. K. Panchal\\
Department of Mathematics, \\ Dr. Babasaheb Ambedkar Marathwada University,\\ Aurangabad-431004 (M.S.) India.\\ pvasudeo.dole@gmail.com; drpanchalsk@gmail.com}
\begin{document}

\maketitle{\bf Abstract:}\\
\indent\qquad In this paper we define canonical sine and cosine transform, convolution operations, prove convolution theorems in space of integrable functions on real space. Further, obtain some results require to construct the spaces of integrable Boehmians then extend this canonical sine and canonical cosine transforms to space of integrable Boehmians and obtain their properties.



\section{Introduction}
\indent\quad The most recent generalizations of functions is the theory of Boehmians. The idea of construction of Boehmians was initiated by the concept of regular operators introduced by Boehme \cite{1}. The space of Boehmians contains all regular operators, all distributions and some objects which are neither operators nor distributions. The construction of Boehmians is given by J. Mikusinski and P. Mikusinski. Boehmians is the new class of generalized functions, has opened the door to new area of research in mathematics \cite{2,3,4,5}. Class of all Boehmians is quite general. P. Mikusinski has studied Fourier transform of integrable Boehmians in \cite{6}.  A. I. Zayed define convolution operations for fractional Fourier transform (FrFT) and extend the FrFT to space of integrable Boehmians by using algebraic technique in \cite{10}. Fourier sine and cosine transforms on Boehmian spaces are studied by R. Roopkumar at.el. in \cite{7}. Fourier sine and cosine transforms are spacial cases of the canonical sine and cosine transforms. The Fourier, Laplace and fractional Fourier transforms are the special cases of linear canonical transform (LCT) and has many applications in several areas, like signal processing and optics \cite{BRY, 9}.

 Let $\mathcal{L}^{1}(\mathbb{R})$ be the space of all complex valued absolutely integrable functions on $\mathbb{R}$ with norm $||f||_{1}=\int_{\mathbb{R}}|f(t)|dt < \infty $. We studied linear canonical transform (LCT) for integrable Boehmians and is defined for $f\in \mathcal{L}^{1}(\mathbb{R})$ in \cite{PS1} as,
\begin{align}
F_{A}=\mathcal{L}_{A}[f(t)](s)=\left\{ {\begin{array}{*{20}{l}}{\frac{1}{\sqrt{2\pi ib}}e^{\frac{i}{2}(\frac{d}{b})s^{2}}\int_{-\infty}^{\infty}e^{\frac{-i}{b}st}e^{\frac{i}{2}\frac{a}{b}t^{2}}f(t)dt \qquad\quad for \quad b\neq 0,}\\
{\sqrt{d}\,e^{\frac{i}{2}cds^{2}} f(ds) \qquad\qquad\qquad\qquad\qquad\quad for \quad b=0,}
\end{array}} \right.\label{1}
\end{align}
where $\mathcal{L}_{A}$ is linear canonical operator with $A=(a, b, c, d)$, parameter $a,b,c,d$ are real number satisfying $ad-bc=1$. The inverse transform for linear canonical transform is given by a linear canonical transform having parameter $A^{-1}=(d, -b, -c, a)$.

 Let $f,g\in \mathcal{L}^{1}(\mathbb{R})$, then the convolution $(f*g)$ is define as $(f*g)(x)=\int_{-\infty}^{\infty}f(t)g(x-t)dt$ is exists and $\|(f*g)\|_{1}=\|f\|_{1} \|g\|_{1}$.

 The canonical cosine and canonical sine transforms are obtain by splinting the definition \eqref{1} of linear canonical transform into real and imaginary part.
 \begin{definition} The canonical cosine transform (CCT) of $f\in\mathcal{L}^{1}(\mathbb{R})$ for $b\neq 0$, is define as
\begin{align}
\mathcal{C}_{A}( f(t))(s)=F_{cA}(s)=\frac{1}{\sqrt{2\pi ib}}e^{\frac{i}{2}(\frac{d}{b})s^{2}}\int_{-\infty}^{\infty}cos(\frac{s}{b}t)e^{\frac{i}{2}(\frac{a}{b})t^{2}}f(t)dt.
\end{align}
where $\mathcal{C}_{A}$ is CCT operator.
\end{definition}
\begin{definition}
The canonical sine transform (CST) of $f\in\mathcal{L}^{1}(\mathbb{R})$ for $b\neq 0$, is define as
\begin{align}
\mathcal{S}_{A}(f(t))(s)=F_{sA}(s)=\frac{1}{\sqrt{2\pi ib}}e^{\frac{i}{2}(\frac{d}{b})s^{2}}\int_{-\infty}^{\infty}sin(\frac{s}{b}t)e^{\frac{i}{2}(\frac{a}{b})t^{2}}f(t)dt,
\end{align}
where $\mathcal{S}_{A}$ is CST operator.
\end{definition}
We see that $f\in\mathcal{L}^{1}(\mathbb{R})$ then $\mathcal{C}_{A}( f(t))$ and $\mathcal{S}_{A}( f(t))$ are the members of $\mathcal{L}^{1}(\mathbb{R})$.
\begin{definition}\cite{Pathak} \textbf{(Regular Distributions)} Let $f$ be the locally integrable function, i.e. absolutely integrable on every finite interval on $\mathbb{R}$, then distribution generated by $f$ is called {\it regular distributions}.
\end{definition}
\section{Canonical Cosine Transform}
\indent\quad In this section we obtain some results which are require to construct space of integrable Boehmians. Let $f, g\in \mathcal{L}^{1}(\mathbb{R})$ and the functions  $\widehat{f}(t)= e^{\frac{i}{2}(\frac{a}{b})t^{2}}f(t)$ and $\widehat{g}(|x-t|)= e^{\frac{i}{2}(\frac{a}{b})(x-t)^{2}}g(x-t)$ , $\widehat{f}(x+t)=e^{\frac{i}{2}(\frac{a}{b})(x+t)^{2}}f(x+t)$. We define the convolution $(f\star g)$ as,
\begin{align}
(f\star g)(t)=\int_{-\infty}^{\infty}\frac{e^{-\frac{i}{2}(\frac{a}{b})t^{2}}}{2}\widehat{f}(x)[\widehat{g}(x+t)+\widehat{g}(|x-t|)] dx.
\end{align}
 Clearly the functions $\widehat{f}, \widehat{g}$ are the members of $\mathcal{L}^{1}(\mathbb{R})$.
\begin{lemma}\label{lemma3.2.1}
 Let $f, g\in \mathcal{L}^{1}(\mathbb{R})$, then the convolutions $(f\star g)$ is also in $\mathcal{L}^{1}(\mathbb{R})$.
\end{lemma}
\textbf{Proof:} Let $f, g\in \mathcal{L}^{1}(\mathbb{R})$  we have,
\begin{align*}
&\|(f\star g)\|_{1}\\
&=\int_{\mathbb{R}}|(f\star g)|dt=\int_{\mathbb{R}}\bigg{|}\int_{-\infty}^{\infty}\frac{e^{-\frac{i}{2}(\frac{a}{b})t^{2}}}{2}\widehat{f}(y)[\widehat{g}(x+t)+\widehat{g}(|x-t|)] dx\bigg{|}dt\\
&\leq \frac{1}{2}\int_{\mathbb{R}}\bigg{|}\int_{-\infty}^{\infty}e^{-\frac{i}{2}(\frac{a}{b})t^{2}}e^{\frac{i}{2}(\frac{a}{b})x^{2}}f(x)e^{\frac{i}{2}(\frac{a}{b})(x+t)^{2}}g(x+t)dx\bigg{|}dt\\
&+\frac{1}{2}\int_{\mathbb{R}}\bigg{|}\int_{-\infty}^{\infty}e^{-\frac{i}{2}(\frac{a}{b})t^{2}}e^{\frac{i}{2}(\frac{a}{b})x^{2}}f(x)e^{\frac{i}{2}(\frac{a}{b})(x-t)^{2}}g(x-t) dx\bigg{|}dt\\
&\leq \frac{1}{2}\int_{\mathbb{R}}\int_{-\infty}^{\infty}|g(x+t)f(x)|dx dt+\frac{1}{2}\int_{\mathbb{R}}\int_{-\infty}^{\infty}|g(x-t)f(x)|dx dt\\
&\leq \frac{1}{2}\int_{\mathbb{R}}\|g(x+t)f(x)\|_{1}dt+\frac{1}{2}\int_{\mathbb{R}}\|g(x-t)f(x)\|_{1}dt\\
&\leq \frac{1}{2}\|f\|_{1}\int_{\mathbb{R}}|g(x+t)|dt+\frac{1}{2}\|f(x)\|_{1}\int_{\mathbb{R}}|g(x-t)|dt= \|f\|_{1} \|g\|_{1}.
\end{align*}$\hfill\blacksquare$
\begin{lemma}\label{3.2.3}
The space $(\mathcal{L}^{1}(\mathbb{R}), \star) $  is commutative semigroup.
\end{lemma}
\textbf{Proof:} Let $f, g, h\in \mathcal{L}^{1}(\mathbb{R})$, from \eqref{lemma3.2.1} we have $(f\star g)\in \mathcal{L}^{1}(\mathbb{R})$.\\
(i) Now we prove commutative property
\begin{align*}
(f\star g)(t)&= \int_{-\infty}^{\infty}\frac{e^{-\frac{i}{2}(\frac{a}{b})t^{2}}}{2}\widehat{f}(x)[\widehat{g}(x+t)+\widehat{g}(|x-t|)] dx\\
&= \frac{1}{2}\int_{-\infty}^{\infty}e^{-\frac{i}{2}(\frac{a}{b})t^{2}}e^{\frac{i}{2}(\frac{a}{b})x^{2}}f(x)e^{\frac{i}{2}(\frac{a}{b})(x+t)^{2}}g(x+t)dx\\
&\quad+ \frac{1}{2}\int_{-\infty}^{\infty}e^{-\frac{i}{2}(\frac{a}{b})t^{2}}e^{\frac{i}{2}(\frac{a}{b})x^{2}}f(x)e^{\frac{i}{2}(\frac{a}{b})(x-t)^{2}}g(x-t) dx= I_{1} + I_{2}
\end{align*}
we put $x+t=y$ in $I_{1}$ and $t-x=y$ in $I_{2}$, we get
\begin{align*}
&= \frac{1}{2}\int_{-\infty}^{\infty}e^{-\frac{i}{2}(\frac{a}{b})t^{2}}e^{\frac{i}{2}(\frac{a}{b})(y-t)^{2}}f(y-t)e^{\frac{i}{2}(\frac{a}{b})y^{2}}g(y)dy\\
&\quad+ \frac{1}{2}\int_{-\infty}^{\infty}e^{-\frac{i}{2}(\frac{a}{b})t^{2}}e^{\frac{i}{2}(\frac{a}{b})(y+t)^{2}}f(y+t)e^{\frac{i}{2}(\frac{a}{b})y^{2}}g(y) dy \\
&=\frac{e^{-\frac{i}{2}(\frac{a}{b})t^{2}}}{2}\bigg{(}\int_{-\infty}^{\infty}\widehat{f}(|y-t|)\widehat{g}(y)dy+\int_{-\infty}^{\infty}\widehat{f}(y+t)\widehat{g}(y)] dy\bigg{)}\\
&=\int_{-\infty}^{\infty}\frac{e^{-\frac{i}{2}(\frac{a}{b})t^{2}}}{2}\widehat{g}(y)[\widehat{f}(|y-t|)+\widehat{f}(y+t)] dx= (g\star f)(t).
\end{align*}
(ii) Associativity
 \begin{align*}
&[f\star (g\star h)](x)= \int_{-\infty}^{\infty}\frac{e^{-\frac{i}{2}(\frac{a}{b})x^{2}}}{2}[\widehat{f}(y+x)+\widehat{f}(|y-x|)] \widehat{(g\star h)}(y)dy\\
&= \int_{-\infty}^{\infty}\frac{e^{-\frac{i}{2}(\frac{a}{b})x^{2}}}{2}e^{\frac{i}{2}(\frac{a}{b})y^{2}}[\widehat{f}(y+x)+\widehat{f}(|y-x|)]  (g\star h)(y)dy\\
&= \int_{-\infty}^{\infty}\frac{e^{-\frac{i}{2}(\frac{a}{b})x^{2}}}{2}e^{\frac{i}{2}(\frac{a}{b})y^{2}}[\widehat{f}(y+x)+\widehat{f}(|y-x|)] \\
&\quad\times\bigg{(}\int_{-\infty}^{\infty}\frac{e^{-\frac{i}{2}(\frac{a}{b})y^{2}}}{2}\widehat{h}(z)[\widehat{g}(z+y)+\widehat{g}(|z-y|)] dz \bigg{)}dy\\
&= \int_{-\infty}^{\infty}\frac{e^{-\frac{i}{2}(\frac{a}{b})x^{2}}}{2}\widehat{h}(z) \bigg{(}\int_{-\infty}^{\infty}\frac{e^{-\frac{i}{2}(\frac{a}{b})y^{2}}}{2}e^{\frac{i}{2}(\frac{a}{b})y^{2}}\\
&\quad\times[\widehat{f}(y+x)+\widehat{f}(|y-x|)][\widehat{g}(z+y)+\widehat{g}(|z-y|)] dy \bigg{)}dz\\
&= \int_{-\infty}^{\infty}\frac{e^{-\frac{i}{2}(\frac{a}{b})x^{2}}}{2}\widehat{h}(z)\bigg{(}\int_{-\infty}^{\infty} \frac{e^{-\frac{i}{2}(\frac{a}{b})y^{2}}}{2}e^{\frac{i}{2}(\frac{a}{b})y^{2}}
[\widehat{f}(y+x)+\widehat{f}(|y-x|)]\widehat{g}(z+y)dy\\
&\quad+\int_{-\infty}^{\infty} \frac{e^{-\frac{i}{2}(\frac{a}{b})y^{2}}}{2}e^{\frac{i}{2}(\frac{a}{b})y^{2}}[\widehat{f}(y+x)+\widehat{f}(|y-x|)]\widehat{g}(|z-y|)dy \bigg{)} dz
\end{align*}
Here we put $y=(u-z)$ in first integration and $y=(u+z)$ in second integration, we get,
 \begin{align*}
&= \int_{-\infty}^{\infty}\frac{e^{-\frac{i}{2}(\frac{a}{b})x^{2}}}{2}\widehat{h}(z)\\
&\quad\bigg{(}\int_{-\infty}^{\infty} \frac{e^{-\frac{i}{2}(\frac{a}{b})(u-z)^{2}}}{2}e^{\frac{i}{2}(\frac{a}{b})(u-z)^{2}}[\widehat{f}(u-z+x)+\widehat{f}(|u-z-x|)]\widehat{g}(u) du\\
&\quad+\int_{-\infty}^{\infty} \frac{e^{-\frac{i}{2}(\frac{a}{b})(u+z)^{2}}}{2}e^{\frac{i}{2}(\frac{a}{b})(u+z)^{2}}[\widehat{f}(u+z+x)+\widehat{f}(|u+z-x|)]\widehat{g}(|u|)du \bigg{)} dz\\
&= \int_{-\infty}^{\infty}\frac{e^{-\frac{i}{2}(\frac{a}{b})x^{2}}}{2}\widehat{h}(z)[e^{\frac{i}{2}(\frac{a}{b})(u-z)^{2}}(f\star g)(u-z)+e^{\frac{i}{2}(\frac{a}{b})(u+z)^{2}}(f\star g)(u+z)]dz\\
&= \int_{-\infty}^{\infty}\frac{e^{-\frac{i}{2}(\frac{a}{b})x^{2}}}{2}\widehat{h}(z)[\widehat{(f\star g)}(u-z)+\widehat{(f\star g)}(u+z)]dz,\,for\, u=x\,we\,get,\\
&=[(f\star g)\star h](x). 
\end{align*}$\hfill\blacksquare$
\begin{theorem}\label{thm3.2.1} (Canonical cosine convolution theorem)
 Let $f, g\in \mathcal{L}^{1} (\mathbb{R})$ and $F_{c}(s), G_{c}(s)$ denote the canonical cosine transform of $f(t)$ and $g(t)$ respectively, then
\begin{align}
\mathcal{C}_{A}[(f\star g)(t)](s)=\sqrt{2\pi ib}\,e^{-\frac{i}{2}(\frac{d}{b})s^{2}} F_{c}(s)G_{c}(s).
\end{align}
\end{theorem}
\textbf{Proof:} Let $f, g\in \mathcal{L}^{1} (\mathbb{R})$ then we have,
\begin{align*}
&\mathcal{C}_{A}[(f\star g)(t)](s)=\frac{e^{\frac{i}{2}(\frac{d}{b})s^{2}}}{\sqrt{2\pi ib}}\int_{-\infty}^{\infty}cos(\frac{s}{b}t)e^{\frac{i}{2}(\frac{a}{b})t^{2}}\\
&\times\bigg{(}\int_{-\infty}^{\infty}\frac{e^{-\frac{i}{2}(\frac{a}{b})t^{2}}}{2}\widehat{f}(x)[\widehat{g}(x+t)+\widehat{g}(|x-t|)] dx \bigg{)}dt\\
&=\frac{e^{\frac{i}{2}(\frac{d}{b})s^{2}}}{2\sqrt{2\pi ib}}\int_{-\infty}^{\infty}e^{\frac{i}{2}(\frac{a}{b})x^{2}}f(x)\int_{-\infty}^{\infty}\\
&\quad\times[e^{\frac{i}{2}(\frac{a}{b})(x+t)^{2}}g(x+t)+ e^{\frac{i}{2}(\frac{a}{b})(t-x)^{2}}g(t-x)]cos(\frac{s}{b}t)dt dx\\
&=\frac{e^{\frac{i}{2}(\frac{d}{b})s^{2}}}{2\sqrt{2\pi ib}}\int_{-\infty}^{\infty}e^{\frac{i}{2}(\frac{a}{b})x^{2}}f(x)\int_{-\infty}^{\infty}e^{\frac{i}{2}(\frac{a}{b})(x+t)^{2}}g(x+t)cos(\frac{s}{b}t)dt dx\\
&\quad +\frac{e^{\frac{i}{2}(\frac{d}{b})s^{2}}}{2\sqrt{2\pi ib}}\int_{-\infty}^{\infty}e^{\frac{i}{2}(\frac{a}{b})x^{2}}f(x)\int_{-\infty}^{\infty}e^{\frac{i}{2}(\frac{a}{b})(t-x)^{2}}g(t-x)cos(\frac{s}{b}t)dt dx= I_{1}+I_{2},
\end{align*}
we put $x+t=y$ in $I_{1}$ and $t-x=y$ in $I_{2}$, we get
\begin{align*}
&=\frac{e^{\frac{i}{2}(\frac{d}{b})s^{2}}}{2\sqrt{2\pi ib}}\int_{-\infty}^{\infty}e^{\frac{i}{2}(\frac{a}{b})x^{2}}f(x)\int_{-\infty}^{\infty}e^{\frac{i}{2}(\frac{a}{b})(y)^{2}}g(y)cos(\frac{s}{b}(y-x))dy dx\\
&\quad +\frac{e^{\frac{i}{2}(\frac{d}{b})s^{2}}}{2\sqrt{2\pi ib}}\int_{-\infty}^{\infty}e^{\frac{i}{2}(\frac{a}{b})x^{2}}f(x)\int_{-\infty}^{\infty}e^{\frac{i}{2}(\frac{a}{b})(y)^{2}}g(y)cos(\frac{s}{b}(y+x))dy dx\\
&=\frac{e^{\frac{i}{2}(\frac{d}{b})s^{2}}}{2\sqrt{2\pi ib}}\int_{-\infty}^{\infty}e^{\frac{i}{2}(\frac{a}{b})x^{2}}f(x)\int_{-\infty}^{\infty}e^{\frac{i}{2}(\frac{a}{b})(y)^{2}}[cos(\frac{s}{b}(y+x))+cos(\frac{s}{b}(y-x))]g(y)dy dx\\
&=\int_{-\infty}^{\infty}e^{\frac{i}{2}(\frac{a}{b})x^{2}}f(x)\bigg{[}\frac{e^{\frac{i}{2}(\frac{d}{b})s^{2}}}{2\sqrt{2\pi ib}}\int_{-\infty}^{\infty}e^{\frac{i}{2}(\frac{a}{b})(y)^{2}}[2cos(\frac{s}{b}y)cos(\frac{s}{b}x)]g(y)dy\bigg{]} dx\\
&=\frac{\sqrt{2\pi ib}}{e^{\frac{i}{2}(\frac{d}{b})s^{2}}}\frac{e^{\frac{i}{2}(\frac{d}{b})s^{2}}}{\sqrt{2\pi ib}}\int_{-\infty}^{\infty}e^{\frac{i}{2}(\frac{a}{b})x^{2}}cos(\frac{s}{b}x)f(x)G_{c}(s)dx\\
&=\sqrt{2\pi ib}\,e^{-\frac{i}{2}(\frac{d}{b})s^{2}}F_{c}(s)G_{c}(s).
\end{align*}$\hfill \blacksquare$

 Let $\bigtriangleup$ be the collection of all sequences of continuous real functions $( \delta_{n} )$ from $\mathcal{L}^{1}(\mathbb{R})$ with the following conditions:\\
 (i) $\int_{\mathbb{R}}e^{(\frac{ia}{b})tx}\delta_{n}(t)dt=1$,  $\forall \, n\in \mathbb{N}$ and $x\in\mathbb{R}$,\\
 (ii) $\|\delta_{n}\|_{1} \leq M$,  $\forall \, n\in \mathbb{N}$ for some $M>0$,\\
 (iii) $\lim_{n\rightarrow \infty}\int_{|t|>\epsilon} |\delta_{n}(t)| dt = 0$, for each $\epsilon >0$.\\
 The members of $\bigtriangleup$ are called {\it delta sequences}.
\begin{example}
The sequence of functions for $a, b, x\in\mathbb{R}; b\neq 0$
\begin{align*}
\delta_{n}(t)=\left\{ {\begin{array}{*{20}{l}}{e^{(\frac{-ia}{b})tx}n^{2}t \qquad\qquad\qquad for \quad 0\leq t\leq \frac{1}{n},}\\
{e^{(\frac{-ia}{b})tx}n^{2} (\frac{2}{n}-t)\qquad\quad\,\, for \quad \frac{1}{n}\leq t\leq \frac{2}{n},}\\
{0 \qquad\qquad\qquad\qquad\qquad\quad otherwise}.
\end{array}} \right.
\end{align*}
\end{example}
\textbf{Solution:} We have to show that the $ ^{1}\delta_{n}(t)$ is member of $\bigtriangleup$.\\
(i) $\int_{\mathbb{R}}e^{(\frac{ia}{b})tx}e^{(\frac{-ia}{b})tx}n^{2}tdt + \int_{\mathbb{R}}e^{(\frac{ia}{b})tx}e^{(\frac{-ia}{b})tx} (\frac{2}{n}-t)dt=\int_{0}^{\frac{1}{n}}n^{2} tdt + \int_{\frac{1}{n}}^{\frac{2}{n}}n^{2} (\frac{2}{n}-t)dt=1$.\\
(ii) $\|\delta_{n}\|_{1}=\int_{0}^{\frac{1}{n}}|2n^{2} t|dt + \int_{\frac{1}{n}}^{\frac{2}{n}}|n^{2} (\frac{2}{n}-t)|dt \leq M$,  sufficiently large $M$ for $\forall \, n\in \mathbb{N}$.\\
(iii) Choose $\epsilon > \frac{2}{n}$, then $ \lim_{n\rightarrow \infty}\int_{0}^{\frac{1}{n}}|n^{2} t| dt +\lim_{n\rightarrow \infty}\int_{\frac{1}{n}}^{\frac{2}{n}}|n^{2} (\frac{2}{n}-t)dt| = 0 $. $ \hfill\blacksquare $

\begin{lemma}\label{lemma3.2.3}
  Let $( \varphi_{n} ), ( \psi_{n} )\in \bigtriangleup$ then $(\varphi_{n}\star\psi_{n}) \in \bigtriangleup$ for each $n\in \mathbb{N}$.
\end{lemma}
\textbf{Proof:} Let $( \varphi_{n} ), ( \psi_{n} )\in \bigtriangleup$ then we have\\
(i)
\begin{align*}
&\int_{\mathbb{R}}e^{(\frac{ia}{b})tx}(\varphi_{n}\star\psi_{n})(t)dt\\
&=\int_{\mathbb{R}}e^{(\frac{ia}{b})tx}\int_{-\infty}^{\infty}\frac{e^{-\frac{i}{2}(\frac{a}{b})t^{2}}}{2}\widehat{\varphi_{n}}(x)[\widehat{\psi_{n}}(t+x)+\widehat{\psi_{n}}(|t-x|)] dx dt\\
&= \frac{1}{2}\int_{\mathbb{R}}e^{(\frac{ia}{b})tx}\int_{-\infty}^{\infty}e^{-\frac{i}{2}(\frac{a}{b})t^{2}}e^{\frac{i}{2}(\frac{a}{b})x^{2}}\varphi_{n}(x)e^{\frac{i}{2}(\frac{a}{b})(t+x)^{2}}\psi_{n}(t+x)dx dt\\
&\quad+\frac{1}{2} \int_{\mathbb{R}}e^{(\frac{ia}{b})tx}\int_{-\infty}^{\infty}e^{-\frac{i}{2}(\frac{a}{b})t^{2}}e^{\frac{i}{2}(\frac{a}{b})x^{2}}\varphi_{n}(x)e^{\frac{i}{2}(\frac{a}{b})(t-x)^{2}}\psi_{n}(t-x) dxdt\\
&= \frac{1}{2}\int_{\mathbb{R}}e^{(\frac{ia}{b})tx}\int_{-\infty}^{\infty}e^{(\frac{ia}{b})x^{2}}\varphi_{n}(x)[e^{(\frac{ia}{b})tx}\psi_{n}(t+x)+e^{(\frac{-ia}{b})tx}\psi_{n}(t-x)] dxdt\\
&= \frac{1}{2}\int_{-\infty}^{\infty}e^{(\frac{ia}{b})tx}\varphi_{n}(x)e^{(\frac{ia}{b})x^{2}}\\
&\quad\times\bigg{[}\int_{\mathbb{R}}e^{(\frac{ia}{b})tx}\psi_{n}(t+x)dt + \int_{\mathbb{R}}e^{(\frac{-ia}{b})tx}\psi_{n}(t-x) dt\bigg{]} dx
\end{align*}
\begin{align*}
&= \frac{1}{2}\int_{-\infty}^{\infty}e^{(\frac{ia}{b})tx}\varphi_{n}(x)e^{(\frac{ia}{b})x^{2}}\\
&\quad\times\bigg{[}\int_{\mathbb{R}}e^{(\frac{ia}{b})x(y-x)}\psi_{n}(y)dy + \int_{\mathbb{R}}e^{(\frac{-ia}{b})x(y+x)}\psi_{n}(y) dy\bigg{]} dx\\
&= \frac{1}{2}\int_{-\infty}^{\infty}e^{(\frac{ia}{b})tx}\varphi_{n}(x)\bigg{[}\int_{\mathbb{R}}e^{(\frac{ia}{b})xy}\psi_{n}(y)dy +\int_{\mathbb{R}}e^{(\frac{-ia}{b})xy}\psi_{n}(y) dy\bigg{]} dx\\
&= \frac{1}{2}\int_{-\infty}^{\infty}e^{(\frac{ia}{b})tx}\varphi_{n}(x)\bigg{[}\int_{\mathbb{R}}e^{(\frac{ia}{b})xy}\psi_{n}(y)dy+\int_{\mathbb{R}}e^{(\frac{ia}{b})xz}\psi_{n}(z) dz\bigg{]} dx=1.\\
\end{align*}
(ii) From lemma \eqref{lemma3.2.1}, we have
\begin{align*}
\|(\varphi_{n}\star\psi_{n})\|_{1}&=\int_{\mathbb{R}}|(\varphi_{n}\star\psi_{n})(t)|dt
\leq \|\varphi_{n}\|_{1}\|\psi_{n}\|_{1}\leq M \quad for \,\, M > 0 .
\end{align*}
(iii)
\begin{align*}
&\lim_{n\rightarrow \infty}\int_{|t|>\epsilon} |(\varphi_{n}\star\psi_{n})(t)| dt \\
&\leq \frac{1}{2}\lim_{n\rightarrow\infty}\int_{|t|>\epsilon}\bigg{(}\int_{-\infty}^{\infty}\big{|}\varphi_{n}(x)\psi_{n}(t+x)\big{|}dxdt+ \int_{-\infty}^{\infty}\big{|}\varphi_{n}(x)\psi_{n}(t-x)\big{|}dxdt\bigg{)}\\
&\leq \lim_{n\rightarrow \infty}\int_{|t|>\epsilon}\|\varphi_{n}(x)\psi_{n}(t+x)\|_{1}dt+\frac{1}{2}\lim_{n\rightarrow \infty}\int_{|t|>\epsilon}\|\varphi_{n}(x)\psi_{n}(t-x)\|_{1}dt\\
&\leq \frac{\|\varphi_{n}\|_{1}}{2}\lim_{n\rightarrow \infty}\int_{|t|>\epsilon}|\psi_{n}(t+x)\big{|}dt+\frac{\|\varphi_{n}\|_{1}}{2}\lim_{n\rightarrow \infty}\int_{|t|>\epsilon}|\psi_{n}(t-x)|dt\\
&\leq \|\varphi_{n}\|_{1}\lim_{n\rightarrow \infty}\int_{|y| > \delta}|\psi_{n}(y)\big{|}dy\longrightarrow 0 \quad for\, each\quad \delta > 0.
\end{align*}
\begin{lemma}\label{lemma3.2.5}
   Let $f \in \mathcal{L}^{1}(\mathbb{R})$ and $( \varphi_{n} )\in \bigtriangleup$  then $(f\star\varphi_{n})\rightarrow f$ as $n\rightarrow \infty$ in $\mathcal{L}^{1}(\mathbb{R})$.
\end{lemma}
\textbf{Proof:} For $\epsilon > 0 $ and $g(x)=f(|x|), \forall x\in\mathbb{R}$ then $g \in \mathcal{L}^{1}(\mathbb{R})$. Since the mapping $t\longmapsto g_{t}$ is uniformly continuous on $\mathbb{R}$, where $g_{t}=g(x-t), \forall x\in\mathbb{R}$, there exist $\delta >0$ such that $\|g_{t}-g_{0}\|_{1} < \epsilon, \forall\,\, t\in [0, \delta)$ (see  Theorem 9.5 \cite{8}). In other words,$\forall\,\, t\in [0, \delta)$, we have
\begin{align*}
\int_{\mathbb{R}}|f(x+t)-f(x)|dx=\int_{\mathbb{R}}|f(|x-t|)-f(|x|)|dx=\|g_{t}-g_{0}\|_{1} < \epsilon
\end{align*}
and  $\int_{\mathbb{R}}[e^{(\frac{ia}{b})x^{2}}f(x+t)-f(x)]dx\leq \int_{\mathbb{R}}[e^{(\frac{ia}{b})x^{2}}f(x+t)-e^{(\frac{ia}{b})x^{2}}f(x)]dx $. Let $f \in \mathcal{L}^{1}(\mathbb{R})$ and $( \varphi_{n} ) \in\bigtriangleup$. Form property (iv) of $\bigtriangleup$, $\varphi_{n}(x)$ has compact support such that $\int_{K}|\varphi_{n}(x)|dx\rightarrow 0\,\, as\,\, n\rightarrow\infty $, Since Supp$\varphi_{n}\rightarrow 0\,\, as\,\, n\rightarrow\infty $, inequality \eqref{3.2.2}, and property (i) of $\bigtriangleup$ we get,
\begin{align*}
&\|(f\star\varphi_{n})-f\|_{1}=\|(\varphi_{n}\star f)-f\|_{1}\\
&=\int_{\mathbb{R}}\Big{|}\int_{-\infty}^{\infty}\frac{e^{-\frac{i}{2}(\frac{a}{b})t^{2}}}{2}\widehat{\varphi_{n}}(x)[\widehat{f}(x+t)+\widehat{f}(|x-t|)]dx-f(t)\Big{|}dt\\
&= \int_{\mathbb{R}}\Bigg{|}\frac{1}{2}\int_{-\infty}^{\infty}e^{-\frac{i}{2}(\frac{a}{b})t^{2}}e^{\frac{i}{2}(\frac{a}{b})x^{2}}\varphi_{n}(x)\\
&\quad [e^{\frac{i}{2}(\frac{a}{b})(x+t)^{2}}f(x+t)+e^{\frac{i}{2}(\frac{a}{b})(|x-t|)^{2}}f(|x-t|)]dx-f(t)\Bigg{|}dt\\
&\leq \int_{\mathbb{R}}\Bigg{|}\frac{1}{2}\int_{-\infty}^{\infty}e^{(\frac{ia}{b})x^{2}}\varphi_{n}(x)[e^{(\frac{ia}{b})xt}f(x+t)+e^{(\frac{ia}{b})xt}f(|x-t|)]dx\\
&\quad-\int_{-\infty}^{\infty}e^{(\frac{ia}{b})xt}\varphi_{n}(x)f(t)dx\Bigg{|}dt
\end{align*}
\begin{align*}
&= \int_{-\infty}^{\infty}\Bigg{|}\frac{1}{2}e^{(\frac{ia}{b})x^{2}}\varphi_{n}(x)\int_{\mathbb{R}}[e^{(\frac{ia}{b})xt}f(x+t)+e^{(\frac{ia}{b})xt}f(|x-t|)]dt\\
&\quad-\varphi_{n}(x)\int_{\mathbb{R}}e^{(\frac{ia}{b})xt}f(t)dt\Bigg{|}dx\\
&\leq \int_{-\infty}^{\infty}\Bigg{|}\frac{1}{2}e^{(\frac{ia}{b})x^{2}}\varphi_{n}(x)\int_{\mathbb{R}}e^{(\frac{ia}{b})xt}[f(x+t)-f(x)+f(|x-t|)-f(t)]dt\Bigg{|}dx\\
&\leq \int_{-\infty}^{\infty}\frac{|\varphi_{n}(x)|}{2}\bigg{(}\int_{\mathbb{R}}|[f(x+t)-f(x)+f(|x-t|)-f(t)]|\bigg{)}dt dx\\
&\leq \int_{-\infty}^{\infty}\frac{|\varphi_{n}(x)|}{2}\bigg{(}\int_{\mathbb{R}}|f(x+t)-f(x)|dt+\int_{\mathbb{R}}|f(|x-t|)-f(t)]|dt\bigg{)} dx\\
&\leq \int_{\mathbb{R}}|\varphi_{n}(x)|\epsilon dx\leq \epsilon\int_{K}|\varphi_{n}(x)|dx\longrightarrow 0\,\, as\,\, n\longrightarrow\infty.
\end{align*}
Hence $(f\star\varphi_{n})\rightarrow f$ as $n\rightarrow \infty$ in $\mathcal{L}^{1}(\mathbb{R})$. $\hfill\blacksquare$
\begin{lemma}\label{lemma3.2.6}
   Let $f_{n}\longrightarrow f$ as $\longrightarrow\infty$ in $\mathcal{L}^{1}(\mathbb{R})$ and $( \varphi_{n} )\in \bigtriangleup$  then $(f_{n}\star\varphi_{n})\rightarrow f$ as $n\rightarrow \infty$ in $\mathcal{L}^{1}(\mathbb{R})$.
\end{lemma}
\textbf{Proof:} For $n\in\mathbb{N}$, we get $\|f_{n}\star \delta_{n}-f\|_{1}=\|f_{n}\star \delta_{n}-f\star \delta_{n}+f\star \delta_{n}-f\|_{1}\leq \|f_{n}\star \delta_{n}-f\star \delta_{n}\|_{1}+\|f\star \delta_{n}-f\|_{1}$. Using lemma \eqref{lemma3.2.4} we get, $\|f_{n}\star \delta_{n}-f\star \delta_{n}\|_{1}\rightarrow 0$ as $n\rightarrow \infty$ and $\|f_{n}\star \delta_{n}-f\|_{1}\leq \|f_{n}\star \delta_{n}-f\star \delta_{n}\|_{1}\leq \|f_{n}-f\|_{1}\|\delta_{n}\|_{1}\leq M\|f_{n}-f\|_{1}\rightarrow 0$ as $n\rightarrow \infty$. where $\|\delta_{n}\|_{1}\leq M$ for $M > 0$.$\hfill\blacksquare$
\begin{lemma}\label{lemma3.2.8}
Let $( \delta_{n} )\in \bigtriangleup$ then $\mathcal{C}_{A}(\delta_{n})$ converges uniformly on each compact set to a constant function $1$ in $\mathcal{L}^{1}(\mathbb{R})$.
\end{lemma}
\textbf{Proof:} Let $( \delta_{n} ) \in \bigtriangleup$. From property (iii) of $\bigtriangleup$ we have  $supp(\delta_{n})\longrightarrow 0$ as $n\longrightarrow \infty$. Let $K$ compact subset of $\mathbb{R}$ such that $\lim_{n\rightarrow \infty}\int_{K}|\delta_{n}(t)|dt\longrightarrow 0$ and also there exist $M>0$ such that $\bigg{|}cos(\frac{s}{b}t)e^{\frac{i}{2}(\frac{a}{b})t^{2}}-\sqrt{2\pi ib}\,e^{\frac{-i}{2}(\frac{d}{b})s^{2}}e^{(\frac{ia}{b})xt}\bigg{|} < M $ for each $|t| > \epsilon > 0$. Hence
\begin{align*}
&\|\big{(}\mathcal{C}_{A}(\delta_{n})-1\big{)}\|_{1}\\
&=\int_{\mathbb{R}}\bigg{|} \frac{e^{\frac{i}{2}(\frac{d}{b})s^{2}}}{\sqrt{2\pi ib}}\int_{-\infty}^{\infty}cos(\frac{s}{b}t)e^{\frac{i}{2}(\frac{a}{b})t^{2}}\delta_{n}(t)dt-\int_{-\infty}^{\infty}e^{(\frac{ia}{b})xt}\delta_{n}(t)dt\bigg{|}ds\\
&=\int_{\mathbb{R}}\bigg{|} \frac{e^{\frac{i}{2}(\frac{d}{b})s^{2}}}{\sqrt{2\pi ib}}\int_{-\infty}^{\infty}\bigg{(}cos(\frac{s}{b}t)e^{\frac{i}{2}(\frac{a}{b})t^{2}}-\sqrt{2\pi ib}\,e^{\frac{-i}{2}(\frac{d}{b})s^{2}}e^{(\frac{ia}{b})xt}\bigg{)}\delta_{n}(t)dt\bigg{|}ds\\
&\leq \frac{M}{\sqrt{2\pi ib}}\int_{\mathbb{R}}\bigg{(}\int_{K}|\delta_{n}(t)|dt\bigg{)} ds \longrightarrow 0 \quad as\,\, n\longrightarrow \infty.
\end{align*}$\hfill\blacksquare$
\section{CCT For Boehmians}
\begin{definition}  A pair of sequences $(( f_{n} ), ( \varphi_{n} ))$ is called a {\it quotient of sequences}, denoted by $f_{n}/\varphi_{n}$, if $( f_{n} )\in \mathcal{L}^{1}(\mathbb{R})$, $( \varphi_{n} )\in\bigtriangleup$ and $f_{m}\star\varphi_{n}=f_{n}\star\varphi_{m}\, \forall\, m, n\in \mathbb{N}$.
\end{definition}
\begin{definition}
Let $( f_{n} ), \{g_{n}\}\in \mathcal{L}^{1}(\mathbb{R})$ and $( \psi_{n} ), ( \varphi_{n} )\in\bigtriangleup$. Two quotient of sequences $f_{n}/ \varphi_{n}$ and $g_{n} / \psi_{n}$ are equivalent if $f_{n}\star\psi_{n}=g_{n}\star\varphi_{n}$ $\forall\, n\in\mathbb{N}$. This is an equivalence relation.
\end{definition}
\begin{definition}
 The equivalence classes of quotient of sequences are called {\it Boehmians}. The space of all Boehmians denoted by $ \mathcal{B}_{\star}=\mathcal{B}(\mathcal{L}^{1}(\mathbb{R}), \bigtriangleup, \star) $ and the members of $\mathcal{B}_{\star}$ is denoted by $F=[f_{n}/\delta_{n}]$. The function $f\in \mathcal{L}^{1}(\mathbb{R})$ then it can be identified with the Boehmian $[f\star\delta_{n}/\delta_{n}]$, where $( \delta_{n} )\in\bigtriangleup$. Let $F=[f_{n}/\delta_{n}]$, then $F\star\delta_{n}=f_{n} \in \mathcal{L}^{1}(\mathbb{R})$ $\forall\, n\in\mathbb{N}$.
\end{definition}
\begin{definition}
A sequence of Boehmians $( f_{n} )$ is called $\Delta-$convergent to a Boehmian $F$ ($\Delta-\lim F_{n}=F$) if there exist a delta sequence $( \delta_{n} )$ such that  $(F_{n}-F)\star\delta_{n}\in \mathcal{L}^{1}(\mathbb{R})$, for every $n\in\mathbb{N}$ and that $\|((F_{n}-F)\star\delta_{n})\|_{1}\rightarrow 0$ as $n\rightarrow \infty$.
\end{definition}
\begin{definition}
 A sequence of Boehmians $( f_{n} )$ is called $\delta-$convergent to a Boehmian $F$ ($\delta-\lim F_{n}=F$) if there exist a delta sequence $( \delta_{n} )$ such that  $F_{n}\star\delta_{k}\in \mathcal{L}^{1}(\mathbb{R})$ and $F\star\delta_{k}\in \mathcal{L}^{1}(\mathbb{R})$ for every $n, k\in\mathbb{N}$ and that $\|((F_{n}-F)\star\delta_{k})\|_{1}\rightarrow 0$ as $n\rightarrow \infty$ for each $k\in\mathbb{N}$.
 \end{definition}
\indent\quad If $( \delta_{n} )$ is a delta sequence, then $[\delta_{n}/\delta_{n}]$ represents an Boehmian for each $n\in\mathbb{N}$. Since the Boehmian [$\delta_{n}/\delta_{n}$] correspond to Dirac delta distribution, all the derivative of $\delta$ are also Boehmian. Further, if $(\delta_{n})$ is infinitely differentiable and bounded support, then the $k^{th}$ derivative of $\delta$ is define as $\delta^{(k)}=[\delta_{n}^{(k)}/\delta_{n}]$, implies  $\delta^{(k)}\in \mathcal{B}_{\star}$, where $k\in\mathbb{N}$. The $k^{th}$ derivative of Boehmian $F\in \mathcal{B}_{\star}$ is define by $F^{(k)}=F\star\delta^{(k)}$. The scalar multiplication, addition and convolution in $\mathcal{B}_{\star}$ is define as below,
\begin{align*}
\lambda[f_{n}/\varphi_{n}]&=[\lambda f_{n}/\varphi_{n}]\\
[f_{n}/\varphi_{n}]+[g_{n}/\psi_{n}]&=[(f_{n}\star\psi_{n}+g_{n}\star\varphi_{n})/\varphi_{n}\star\psi_{n}]\\
[f_{n}/\varphi_{n}]\star[g_{n}/\psi_{n}]&=[f_{n}\star g_{n}/\varphi_{n}\star\psi_{n}].
\end{align*}

Let $\bigtriangleup_{0}=\{\mathcal{C}_{A}(\delta_{n}); ( \delta_{n} )\in \bigtriangleup \}$ be the space of complex valued functions on $(\mathbb{R})$, the operation $\cdot$ is pointwise multiplication and $C_{0}(\mathbb{R})$ be the space of all continuous functions vanishing at infinity on $(\mathbb{R})$. Now we construct the another space of Boehmians, denoted by $\mathcal{B}= (\mathcal{L}^{1}(\mathbb{R}), C_{0}(\mathbb{R})\cap \mathcal{L}^{1}(\mathbb{R}), \bigtriangleup_{0}, \cdot)$. This is the range of CCT on $\mathcal{B}_{\star}$ and each element of $\mathcal{B}$ is denoted by $\mathcal{C}_{A} (f_{n})/ \mathcal{C}_{A} (\delta_{n})$ for all $n\in\mathbb{N}$, where $( f_{n} )\in \mathcal{L}^{1}(\mathbb{R}), ( \delta_{n} )\in \bigtriangleup$.

\begin{lemma} \label{lemma3.3.1} Let $f, g\in \mathcal{L}^{1}(\mathbb{R})$ and $ \varphi, \psi \in C_{0}(\mathbb{R})$ then for any $\lambda\in\mathbb{C}$,\\ (i) $f\cdot\varphi \in \mathcal{L}^{1}(\mathbb{R})$;\\ (ii) $(f+g)\cdot\varphi =f\cdot\varphi + f\cdot\varphi$;\\ (iii) $(\lambda f)\cdot\varphi=\alpha(f\cdot\varphi)$\\ (iv) $f\cdot(\varphi\cdot\psi)=(f\cdot\varphi)\cdot\psi$.
\end{lemma}
\textbf{Proof:} Let $f, g\in \mathcal{L}^{1}(\mathbb{R}); \varphi, \psi \in C_{0}(\mathbb{R})$ and $\lambda\in\mathbb{C}$, we have\\
(i) $\|f\cdot\varphi\|_{1}=\int_{-\infty}^{\infty}|f(x)\varphi(x)|dx=Max_{x\in\mathbb{R}}|\varphi(x)|\int_{-\infty}^{\infty}|f(x)|^{2}dx < \infty$, since $\varphi(x)$ vanishes at infinity. This implies $f\cdot\varphi \in \mathcal{L}^{1}(\mathbb{R})$.\\
(ii) $(f(x)+g(x))\cdot\varphi(x)=f(x)\cdot\varphi(x) + g(x)\cdot\varphi(x)$, since pointwise multiplication is distributive over addition.\\
(iii) $(\lambda f(x))\cdot\varphi(x)=\lambda(f(x)\cdot\varphi(x))$, since scalar multiplication is associative.\\
(iv) $f\cdot(\varphi\cdot\psi)=(f\cdot\varphi)\cdot\psi$, since pointwise multiplication is associative.$\hfill\blacksquare$
\begin{lemma} \label{lemma3.3.2} Let $f_{n}\rightarrow f$ as $n\rightarrow\infty$ in $\mathcal{L}^{1}(\mathbb{R})$ and $\varphi\in C_{0}(\mathbb{R})$ then $f_{n}\cdot\varphi\rightarrow f\cdot\varphi$ as $n\rightarrow\infty$ in $\mathcal{L}^{1}(\mathbb{R})$.
\end{lemma}
\textbf{Proof:} Let $f_{n}\rightarrow f$ as $n\rightarrow\infty$ in $\mathcal{L}^{1}(\mathbb{R})$ and $\varphi\in C_{0}(\mathbb{R})$, we have $\|f_{n}\cdot\varphi-f\cdot\varphi\|_{1}=\|(f_{n}-f)\cdot\varphi\|_{1}\leq Max_{x\in\mathbb{R}}|\varphi(x)|\|(f_{n}-f)\|_{1}\rightarrow 0$ as $ n\rightarrow\infty $.$\hfill\blacksquare$
\begin{lemma} Let $f_{n}\longrightarrow f$ as $n\longrightarrow\infty$ in $\mathcal{L}^{1}(\mathbb{R})$ and $\mathcal{C}_{A}(\varphi_{n})\in \bigtriangleup_{0}$ then $f_{n}\cdot \mathcal{C}_{A}(\varphi_{n})\rightarrow f$ in $\mathcal{L}^{1}(\mathbb{R})$.
\end{lemma}
\textbf{Proof:} Let $f_{n}\rightarrow f$ as $n\longrightarrow\infty$ in $\mathcal{L}^{1}(\mathbb{R})$ and $\mathcal{C}_{A}(\varphi_{n})\in \bigtriangleup_{0}$. Using lemma \eqref{lemma3.2.8} we get,
 \begin{align*}
&\|f_{n}\cdot \mathcal{C}_{A}(\varphi_{n})- f\|_{1}=\|f_{n}\cdot \mathcal{C}_{A}(\varphi_{n})- f\cdot \mathcal{C}_{A}(\varphi_{n})+ f\cdot \mathcal{C}_{A}(\varphi_{n})-f\|_{1}\\
&\leq \|f_{n}-f\|_{1}\|\mathcal{C}_{A}(\varphi_{n})\|_{1}+ \|f\|_{1} \|\mathcal{C}_{A}(\varphi_{n})-1\|_{1}\longrightarrow 0\quad as\,\,n\longrightarrow \infty.
 \end{align*}$\hfill\blacksquare$
\begin{lemma} \label{lemma3.3.4} Let $ \mathcal{C}_{A}(\varphi_{n}), \mathcal{C}_{A}(\psi_{n}) \in \bigtriangleup_{0}$ then $\mathcal{C}_{A}(\varphi_{n})\cdot \mathcal{C}_{A}(\psi_{n})\in \bigtriangleup_{0}$
\end{lemma}
\textbf{Proof:} Let $ \mathcal{C}_{A}(\varphi_{n}), \mathcal{C}_{A}(\psi_{n}) \in \bigtriangleup_{0}$. From the theorem \eqref{thm3.2.1} and lemma \eqref{lemma3.2.3} we get $\mathcal{C}_{A}(\varphi_{n})\cdot \mathcal{C}_{A}(\psi_{n})= \frac{e^{\frac{i}{2}(\frac{d}{b})s^{2}}}{\sqrt{2\pi ib}} \mathcal{C}_{A}(\varphi_{n}\star\psi_{n})\in \bigtriangleup_{0}$. $\hfill\blacksquare$
\begin{definition}
Let $( f_{n} )\in \mathcal{L}^{1}(\mathbb{R})$ and $( \delta_{n} )\in\bigtriangleup$ for $n\in\mathbb{N}$, we define the canonical cosine transform $\mathcal{C}_{A} : \mathcal{B}_{\star}\longrightarrow\mathcal{B}$ as
 \begin{align}
 \mathcal{C}_{A} [f_{n}/\delta_{n}]= \mathcal{C}_{A} (f_{n})/ \mathcal{C}_{A} (\delta_{n}) \qquad for \quad [f_{n}/\delta_{n}]\in \mathcal{B}_{\star}.
\end{align}
\end{definition}
\indent\quad The CCT on $\mathcal{B}_{\star}$ is well-defined. Indeed if  $[f_{n}/\delta_{n}]\in \mathcal{B}_{\star}$, then $f_{n}\star\delta_{m}=f_{m}\star\delta_{n}$ for all $m, n \in \mathbb{N}$. Applying the CCT on both sides, we get $\mathcal{C}_{A} (f_{n}) \mathcal{C}_{A} (\delta_{m})=\mathcal{C}_{A} (f_{m}) \mathcal{C}_{A} (\delta_{n})$ for all $m, n \in \mathbb{N}$ and hence $\mathcal{C}_{A} (f_{n})/ \mathcal{C}_{A} (\delta_{n})\in \mathcal{B}$. Further if $[f_{n}/\psi_{n}]=[g_{n}/\delta_{n}]\in \mathcal{B}_{\star}$ then we have $f_{n}\star\delta_{n}=g_{n}\star\psi_{n}$ for all $ n \in \mathbb{N}$. Again applying the CCT on both sides, we get $\mathcal{C}_{A} (f_{n}) \mathcal{C}_{A} (\delta_{n})=\mathcal{C}_{A} (g_{n}) \mathcal{C}_{A} (\psi_{n})$ for all $ n \in \mathbb{N}$. i.e. $\mathcal{C}_{A} (f_{n})/ \mathcal{C}_{A} (\psi_{n})=\mathcal{C}_{A} (g_{n})/ \mathcal{C}_{A} (\delta_{n})$ in $\mathcal{B}$.
\begin{lemma}\label{lemma3.5}
Let $( f_{n} )\in \mathcal{L}^{1}(\mathbb{R})$ for $n\in\mathbb{N}$, then
\begin{align}
\mathcal{C}_{A} [f_{n}](s)= \sqrt{\frac{1}{2\pi ib}}e^{\frac{i}{2}(\frac{d}{b})s^{2}}\int_{-\infty}^{\infty}cos(\frac{s}{b}t)e^{\frac{i}{2}(\frac{a}{b})t^{2}}f_{n}(t)dt
\end{align}
converges uniformly on each compact subset of $\mathbb{R}$.
\end{lemma}
\textbf{Proof:} Let $( \delta_{n} )$ is delta sequence, then $\mathcal{C}_{A}[\delta_{n}]$ converges uniformly on each compact set to a function $1$. Hence, for each compact subset $M$, $\mathcal{C}_{A}[\delta_{m}] > 0$  on $M$, for almost all $m\in\mathbb{N}$ Moreover,
\begin{align}
\mathcal{C}_{A} (f_{n})&= \mathcal{C}_{A} (f_{n})\frac{\mathcal{C}_{A} (\delta_{m})}{\mathcal{C}_{A} (\delta_{m})}=\frac{e^{\frac{i}{2}(\frac{d}{b})s^{2}}}{\sqrt{2\pi ib}}\frac{\mathcal{C}_{A}(f_{n}\star\delta_{m})}{\mathcal{C}_{A} (\delta_{m})}\nonumber\\
&=\frac{e^{\frac{i}{2}(\frac{d}{b})s^{2}}}{\sqrt{2\pi ib}}\frac{\mathcal{C}_{A}(f_{m}\star\delta_{n})}{\mathcal{C}_{A} (\delta_{m})}=\frac{\mathcal{C}_{A}(f_{m})}{\mathcal{C}_{A}(\delta_{m})}{\mathcal{C}_{A} (\delta_{n})},\qquad on\quad M,
\end{align}
as $n\rightarrow \infty$ we get $\mathcal{C}_{A} (f_{n})\rightarrow \frac{\mathcal{C}_{A}(f_{m})}{\mathcal{C}_{A}(\delta_{m})}= F $, on each compact subset $M$ of $\mathbb{R}$, where $[\frac{f_{m}}{\delta_{m}}]=F$ for each $m\in\mathbb{N}$. $\hfill\blacksquare$

In view of the above lemma, the canonical cosine transform of Boehmian $F\in \mathcal{B}_{\star}$ is define as
\begin{align}
F=\lim_{n\rightarrow\infty}\mathcal{C}_{A} (f_{n}).
\end{align}
Now we have to show that the above definition is well defined. Let two Boehmians $F=[f_{n}/ \varphi_{n}]$ and $G=[g_{n} / \psi_{n}]$ are the equivalent. i.e. $f_{n} \star \psi_{n}=g_{n} \star \varphi_{n}$ for all $n\in\mathbb{N}$. Taking canonical cosine transform on both side, we have
\begin{align*}
 \mathcal{C}_{A}[f_{n}\star\psi_{n}](s)&=\mathcal{C}_{A}[g_{n}\star\varphi_{n}](s)\\
\mathcal{C}_{A} [f_{n}](s)\mathcal{C}_{A}[\psi_{n}](s)&= \mathcal{C}_{A} [g_{n}](s)\mathcal{C}_{A}[\varphi_{n}](s)
\end{align*}
which implies that
\begin{align*}
\lim_{n\rightarrow \infty}\ \mathcal{C}_{A}[f_{n}](s)&=\lim_{n\rightarrow \infty}\mathcal{C}_{A}[g_{n}](s)\\
\mathcal{C}_{A}[F](s)&=\mathcal{C}_{A}[G](s).
\end{align*}
\begin{theorem}
The canonical cosine transform $\mathcal{C}_{A}:\mathcal{B}_{\star}\longrightarrow \mathcal{B}$ is consistent with $\mathcal{C}_{A}: \mathcal{L}^{1}(\mathbb{R})\longrightarrow \mathcal{L}^{1}(\mathbb{R})$.
\end{theorem}
\textbf{Proof:} Let $f\in \mathcal{L}^{1}(\mathbb{R})$. The Boehmian representing $f$ in $\mathcal{B}_{\star}$ is $[(f\star\delta_{n})/\delta_{n}]$ where $( \delta_{n} )\in \bigtriangleup$. Then for each $n\in\mathbb{N}$, $\mathcal{C}_{A}((f\star\delta_{n})/\delta_{n})= \frac{\sqrt{2\pi ib}}{e^{\frac{i}{2}(\frac{d}{b})s^{2}}}\mathcal{C}_{A} (f) \mathcal{C}_{A} (\delta_{n})/\mathcal{C}_{A} (\delta_{n})$, which is a Boehmian in $\mathcal{B}$, which represents $\frac{\sqrt{2\pi ib}}{e^{\frac{i}{2}(\frac{d}{b})s^{2}}}\mathcal{C}_{A}(f)$.$\hfill\blacksquare$
\begin{theorem}
The canonical cosine transform $\mathcal{C}_{A}:\mathcal{B}_{\star}\longrightarrow \mathcal{B}$ is a bijection.
\end{theorem}
\textbf{Proof:} Let $F=[f_{n}/\varphi_{n}], G=[g_{m}/\psi_{m}]\in \mathcal{B}_{\star}$ such that $\mathcal{C}_{A}(F)=\mathcal{C}_{A}(G)$. From this we get $\mathcal{C}_{A}(f_{n}) \mathcal{C}_{A}(\psi_{m})=\mathcal{C}_{A}(g_{m}) \mathcal{C}_{A}(\varphi_{n})$ for all $m, n\in\mathbb{N}$ and hence $\mathcal{C}_{A}(f_{n}\star\psi_{m})=\mathcal{C}_{A}(g_{m}\star\varphi_{n})$ for all $m, n\in\mathbb{N}$. Since CCT is one-to-one in $\mathcal{L}^{1}(\mathbb{R})$, we get
$(f_{n}\star\psi_{m})=(g_{m}\star\varphi_{n})$ for all $m, n\in\mathbb{N}$. This implies $F=G$.
\par Let $G=[g_{m}/\psi_{m}]\in \mathcal{B}_{\star}$. Since $\mathcal{C}_{A}: \mathcal{L}^{1}(\mathbb{R})\longrightarrow \mathcal{L}^{1}(\mathbb{R})$ is onto. Choose $f_{n}\in \mathcal{L}^{1}(\mathbb{R})$ such that $g_{n}= \mathcal{C}_{A} (f_{n})$. Now using the relation $g_{n}\cdot \mathcal{C}_{A}(\psi_{m})=g_{m}\cdot \mathcal{C}_{A}(\psi_{n})$ for all $m, n\in\mathbb{N}$. we obtain $\mathcal{C}_{A}(f_{n})\cdot \mathcal{C}_{A}(\psi_{m})=\mathcal{C}_{A} (f_{m})\cdot \mathcal{C}_{A}(\psi_{n}) \Rightarrow \mathcal{C}_{A}(f_{n}\star\psi_{m})=\mathcal{C}_{A} (f_{m}\star\psi_{n})$ Since CCT is one-to-one in $\mathcal{L}^{1}(\mathbb{R})$, we get $(f_{n}\star\psi_{m})=(f_{m}\star\psi_{n})$. Now if we take $F=[f_{n}/\psi_{n}]$ then $F\in \mathcal{B}_{\star}$ and $\mathcal{C}_{A}(F)=G$. Thus the theorem is hold. $\hfill\blacksquare$
\begin{theorem}
Let $F, G \in \mathcal{B}_{\star}$ then\\
(1) $\mathcal{C}_{A}(\lambda F + G)= \lambda \mathcal{C}_{A} (F)+ \mathcal{C}_{A} (G)$, for any complex $\lambda$,\\
(2) $\mathcal{C}_{A} (F\star G)= \mathcal{C}_{A}(F) \mathcal{C}_{A}(G)$,\\
(3) $\mathcal{C}_{A} (F(kt))=\frac{1}{k}e^{(1-\frac{1}{k^{2}})\frac{i}{2}\frac{d}{b}s^{2}}\mathcal{C}_{A} (e^{(\frac{1}{k^{2}}-1)\frac{i}{2}\frac{a}{b}t^{2}}F)(\frac{s}{k})$\\
(4) $\mathcal{C}_{A} (F (x+\tau))=e^{\frac{i}{2}(\frac{a}{b})\tau^{2}}\{cos(\frac{s}{b}\tau)\mathcal{C}_{A}(e^{-it\tau(\frac{a}{b})} F)+i sin(\frac{s}{b}\tau)\mathcal{S}_{A}(e^{-it\tau(\frac{a}{b})} F)\}$\\
(5) $\mathcal{C}_{A}(e^{ixt} F)(s)=\mathcal{C}_{A}(cos(xt) F)(s)+i \mathcal{C}_{A} (sin(xt) F)(s)$\\
(6) $\mathcal{C}_{A}\big{(}cos (xt)F\big{)}(s)=\frac{e^{\frac{-i}{2}(db)x^{2}}}{2}\big{(}e^{-idxs}\mathcal{C}_{A} \big{(}F\big{)}(s+bx)+e^{idxs}\mathcal{C}_{A} \big{(}F\big{)}(s-bx)\big{)} $\\
(7) $\mathcal{C}_{A}\big{(}sin (xt)F\big{)}=\frac{e^{\frac{-i}{2}(db)x^{2}}}{2}\big{(}e^{-idxs}\mathcal{S}_{A} \big{(}F\big{)}(s+bx)-e^{idxs}\mathcal{S}_{A} \big{(}F\big{)}(s-bx)\big{)} $.
\end{theorem}
\textbf{Proof:} Let $F, G \in \mathcal{B}_{\star}$ and $\lambda $ be any scalar,\\
(1)Linearity
\begin{align*}
\mathcal{C}_{A}[(\lambda F + G)(t)](s)&=\frac{1}{\sqrt{2\pi ib}}e^{\frac{i}{2}(\frac{d}{b})s^{2}}\int_{-\infty}^{\infty}cos(\frac{s}{b}t)e^{\frac{i}{2}(\frac{a}{b})t^{2}}(\lambda F + G)(t)dt\\
&=\frac{1}{\sqrt{2\pi ib}}e^{\frac{i}{2}(\frac{d}{b})s^{2}}\lambda\int_{-\infty}^{\infty}cos(\frac{s}{b}t)e^{\frac{i}{2}(\frac{a}{b})t^{2}} F(t)dt\\
&\quad+\frac{1}{\sqrt{2\pi ib}}e^{\frac{i}{2}(\frac{d}{b})s^{2}}\int_{-\infty}^{\infty}cos(\frac{s}{b}t)e^{\frac{i}{2}(\frac{a}{b})t^{2}}G(t)dt\\
&=\lambda\mathcal{C}_{A}[ F (t)](s)+\mathcal{C}_{A}[ G (t)](s).
\end{align*}
(2) Convolution of Boehmians
\begin{align*}
\mathcal{C}_{A} (F \star G)&=\mathcal{C}_{A} ([f_{n}/\varphi_{n}]\star[g_{n}/\psi_{n}])=\mathcal{C}_{A} [(f_{n}\star g_{n})/(\varphi_{n}\star\psi_{n})]\\
&=\mathcal{C}_{A} [(f_{n}\star g_{n})]/\mathcal{C}_{A} [(\varphi_{n}\star\psi_{n})]\\
&=\frac{\mathcal{C}_{A}(f_{n})\mathcal{C}_{A}(g_{n})}{\mathcal{C}_{A}(\varphi_{n})\mathcal{C}_{A}(\psi_{n})}=\frac{\mathcal{C}_{A}(f_{n})}{\mathcal{C}_{A}(\varphi_{n})}\frac{\mathcal{C}_{A}(g_{n})}{\mathcal{C}_{A}(\psi_{n})}\\
&=\mathcal{C}_{A}[f_{n}/\varphi_{n}]\mathcal{C}_{A}[g_{n}/\psi_{n}]=\mathcal{C}_{A}(F)\mathcal{C}_{A}(G).
\end{align*}
(3) Scaling property,  for $k\in\mathbb{R}$ we have
\begin{align*}
&\mathcal{C}_{A}[F(kt)](u)=\frac{1}{\sqrt{2\pi ib}}e^{\frac{i}{2}(\frac{d}{b})s^{2}}\int_{-\infty}^{\infty}cos(\frac{s}{b}t)e^{\frac{i}{2}\frac{a}{b}t^{2}}F(kt)dt\\
&=\frac{1}{k\sqrt{2\pi ib}}e^{\frac{i}{2}(\frac{d}{b})s^{2}}\int_{-\infty}^{\infty}cos(\frac{s}{kb}x)e^{\frac{i}{2}\frac{a}{b}(\frac{x}{k})^{2}}F(x)dx\\
&=\frac{e^{(1-\frac{1}{k^{2}})\frac{i}{2}\frac{d}{b}s^{2}}}{k\sqrt{2\pi ib}}e^{\frac{i}{2}(\frac{d}{b})(\frac{s}{k})^{2}}\int_{-\infty}^{\infty}cos\bigg{(}\frac{s}{k}\frac{1}{b}x\bigg{)}e^{\frac{i}{2}\frac{a}{b}x^{2}}e^{(\frac{1}{k^{2}}-1)\frac{i}{2}\frac{a}{b}x^{2}}F(x)dx\\
&=\frac{1}{k}e^{(1-\frac{1}{k^{2}})\frac{i}{2}\frac{d}{b}s^{2}}\mathcal{C}_{A} (e^{(\frac{1}{k^{2}}-1)\frac{i}{2}\frac{a}{b}t^{2}}F)(\frac{s}{k}).
\end{align*}
(4) Translation property,  for $\tau\in\mathbb{R}$ we have
\begin{align*}
&\mathcal{C}_{A} [F (t+\tau)](u)=\frac{1}{\sqrt{2\pi ib}}e^{\frac{i}{2}(\frac{d}{b})s^{2}}\int_{-\infty}^{\infty}cos(\frac{s}{b}t)e^{\frac{i}{2}\frac{a}{b}t^{2}}F(t+\tau)dt\\
&=\frac{1}{\sqrt{2\pi ib}}e^{\frac{i}{2}(\frac{d}{b})s^{2}}\int_{-\infty}^{\infty}cos(\frac{s}{b}(x-\tau))e^{\frac{i}{2}\frac{a}{b}(x-\tau)^{2}}F(x)dx\\
&=\frac{e^{\frac{i}{2}(\frac{d}{b})s^{2}}}{\sqrt{2\pi ib}}\int_{-\infty}^{\infty}e^{\frac{i}{2}\frac{a}{b}(x-\tau)^{2}}[cos(\frac{s}{b}x)cos(\frac{s}{b}\tau)+sin(\frac{s}{b}x)sin(\frac{s}{b}\tau)]F(x)dx\\
&=\frac{e^{\frac{i}{2}(\frac{d}{b})s^{2}}}{\sqrt{2\pi ib}}\int_{-\infty}^{\infty}e^{\frac{i}{2}\frac{a}{b}(t-\tau)^{2}}[cos(\frac{s}{b}t)cos(\frac{s}{b}\tau)+sin(\frac{s}{b}t)sin(\frac{s}{b}\tau)]F(t)dt\\
&=e^{\frac{i}{2}(\frac{a}{b})\tau^{2}}\{cos(\frac{s}{b}\tau)\mathcal{C}_{A}(e^{-it\tau(\frac{a}{b})} F)(s)+i sin(\frac{s}{b}\tau)\mathcal{S}_{A}(e^{-it\tau(\frac{a}{b})} F)(s)\}.
\end{align*}
(5)
\begin{align*}
 &\mathcal{C}_{A}(e^{ixt} F)(s)=\frac{1}{\sqrt{2\pi ib}}e^{\frac{i}{2}(\frac{d}{b})s^{2}}\int_{-\infty}^{\infty}cos(\frac{s}{b}t)e^{\frac{i}{2}\frac{a}{b}t^{2}}e^{ixt} F(t)dt\\
 &=\frac{1}{\sqrt{2\pi ib}}e^{\frac{i}{2}(\frac{d}{b})s^{2}}\int_{-\infty}^{\infty}cos(\frac{s}{b}t)e^{\frac{i}{2}\frac{a}{b}t^{2}}[cos(xt)+isin(xt)] F(t)dt\\
 &=\mathcal{C}_{A}(cos(xt) F)(s)+i \mathcal{C}_{A} (sin(xt) F)(s).
\end{align*}
(6)
\begin{align*}
&\mathcal{C}_{A}(cos(xt) F(t))(s)=\frac{1}{\sqrt{2\pi ib}}e^{\frac{i}{2}(\frac{d}{b})s^{2}}\int_{-\infty}^{\infty}cos(\frac{s}{b}t)e^{\frac{i}{2}\frac{a}{b}t^{2}}cos(xt) F(t)dt\\
&=\frac{1}{2\sqrt{2\pi ib}}e^{\frac{i}{2}(\frac{d}{b})s^{2}}\int_{-\infty}^{\infty}e^{\frac{i}{2}\frac{a}{b}t^{2}}[cos(\frac{(s+bx)}{b}t)+cos(\frac{(s-bx)}{b}t)]F(t)dt\\
&=\frac{e^{\frac{i}{2}(\frac{d}{b})(bx)^{2}+2bx}}{2\sqrt{2\pi ib}}e^{\frac{i}{2}(\frac{d}{b})(s+bx)^{2}}\int_{-\infty}^{\infty}e^{\frac{i}{2}\frac{a}{b}t^{2}}cos(\frac{(s+bx)}{b}t)F(t)dt\\
&\quad+\frac{e^{\frac{i}{2}(\frac{d}{b})(bx)^{2}-2bx}}{2\sqrt{2\pi ib}}e^{\frac{i}{2}(\frac{d}{b})(s-bx)^{2}}\int_{-\infty}^{\infty}e^{\frac{i}{2}\frac{a}{b}t^{2}}cos(\frac{(s-bx)}{b}t)F(t)dt\\
&=\frac{e^{\frac{i}{2}(db)x^{2}}}{2}\big{(}e^{idxs}\mathcal{C}_{A} \big{(}F\big{)}(s+bx)+e^{-idxs}\mathcal{C}_{A} \big{(}F\big{)}(s-bx)\big{)}.
\end{align*}
(7)
\begin{align*}
&\mathcal{C}_{A}(sin(xt) F(t))(s)=\frac{1}{\sqrt{2\pi ib}}e^{\frac{i}{2}(\frac{d}{b})s^{2}}\int_{-\infty}^{\infty}cos(\frac{s}{b}t)e^{\frac{i}{2}\frac{a}{b}t^{2}}sin(xt) F(t)dt\\
&=\frac{1}{2\sqrt{2\pi ib}}e^{\frac{i}{2}(\frac{d}{b})s^{2}}\int_{-\infty}^{\infty}e^{\frac{i}{2}\frac{a}{b}t^{2}}[sin(\frac{(s+bx)}{b}t)-sin(\frac{(s-bx)}{b}t)]F(t)dt\\
&=\frac{e^{\frac{i}{2}(\frac{d}{b})(bx)^{2}+2bx}}{2\sqrt{2\pi ib}}e^{\frac{i}{2}(\frac{d}{b})(s+bx)^{2}}\int_{-\infty}^{\infty}e^{\frac{i}{2}\frac{a}{b}t^{2}}sin(\frac{(s+bx)}{b}t)F(t)dt\\
&\quad-\frac{e^{\frac{i}{2}(\frac{d}{b})(bx)^{2}-2bx}}{2\sqrt{2\pi ib}}e^{\frac{i}{2}(\frac{d}{b})(s-bx)^{2}}\int_{-\infty}^{\infty}e^{\frac{i}{2}\frac{a}{b}t^{2}}sin(\frac{(s-bx)}{b}t)F(t)dt\\
&=\frac{e^{\frac{i}{2}(db)x^{2}}}{2}\big{(}e^{idxs}\mathcal{S}_{A} \big{(}F\big{)}(s+bx)-e^{-idxs}\mathcal{S}_{A} \big{(}F\big{)}(s-bx)\big{)}.
\end{align*}
$\hfill\blacksquare$
\begin{theorem}
If $\delta-\lim F_{n}=F$, then $\mathcal{C}_{A} (F_{n})\rightarrow \mathcal{C}_{A} (F)$ as $n\longrightarrow \infty$ on each compact subset of $\mathbb{R}$.
\end{theorem}
\textbf{Proof:} Let $(\delta_{m})$ be a delta sequence such that $F_{n}\star\delta_{m}, F\star\delta_{m}\in \mathcal{L}^{1}(\mathbb{R}) $ for all $n,m\in \mathbb{N}$ and $\|((F_{n}-F)\star\delta_{m})\|_{1}\rightarrow 0$ as $n\rightarrow \infty$ for each $m\in\mathbb{N}$. Let $M$ be a compact subset of $\mathbb{R}$ then $\mathcal{C}_{A}(\delta_{m})> 0$ on $M$ for almost all $m\in\mathbb{N}$. Since $\mathcal{C}_{A}(\delta_{m})$ is a continuous function and $\mathcal{C}_{A}(F_{n})\star \mathcal{C}_{A}(\delta_{m})- \mathcal{C}_{A}(F)\star \mathcal{C}_{A}(\delta_{m})=\big{(}(\mathcal{C}_{A}(F_{n})- \mathcal{C}_{A}(F))\big{)}\star \mathcal{C}_{A}(\delta_{m}))$, implies $\|\big{(}(\mathcal{C}_{A}(F_{n})- \mathcal{C}_{A}(F))\star \mathcal{C}_{A}(\delta_{m})\big{)}\|_{1}\rightarrow 0$ as $n\rightarrow \infty$, thus $\mathcal{C}_{A}(F_{n})\rightarrow \mathcal{C}_{A}(F)$ as $n\rightarrow \infty$ on each $M. \hfill\blacksquare $

\section{Canonical Sine Transform}
\indent\quad In this section we define convolution operation $\Theta$, prove the convolution theorem and Plancherel type theorem for canonical sine transform (CST). Also obtain some result which are require to construct  space of Boehmians. Let $f, g\in \mathcal{L}^{1}(\mathbb{R})$ and the functions  $\widehat{f}(t)= e^{\frac{i}{2}(\frac{a}{b})t^{2}}f(t)$ and $\widehat{g}(|x-t|)= \frac{e^{\frac{i}{2}(\frac{a}{b})(x-t)^{2}}}{2} g(x-t) $ , $\widehat{g}(x+t)= e^{\frac{i}{2}(\frac{a}{b})(x+t)^{2}}g(x+t)$. We define convolution as,
\begin{align}
(f\Theta g)(t)=\int_{-\infty}^{\infty}\frac{e^{-\frac{i}{2}(\frac{a}{b})t^{2}}}{2}\widehat{f}(x)[\widehat{g}(|x-t|)-\widehat{g}(x+t)] dx
\end{align}
 and
\begin{align}
f \Theta g= f\star g- f \otimes g,\label{3.4.1}
\end{align}
where $(f\otimes g)(t) = \int_{\mathbb{R}}e^{-\frac{i}{2}(\frac{a}{b})t^{2}}\widehat{f}(x)\widehat{g}(x+t) dx $ for all $x\in\mathbb{R}$.
\begin{lemma}
The space $(\mathcal{L}^{1}(\mathbb{R}), \otimes) $  is commutative semigroup.
\end{lemma}
\textbf{Proof:}  Let $f, g, h\in \mathcal{L}^{1}(\mathbb{R})$,
\begin{align*}
(f\otimes g)(t) &= \int_{\mathbb{R}}e^{-\frac{i}{2}(\frac{a}{b})t^{2}}\widehat{f}(x)\widehat{g}(x+t) dx \\
&= \int_{\mathbb{R}}e^{-\frac{i}{2}(\frac{a}{b})t^{2}}\widehat{f}(x-t)\widehat{g}(x) dx\qquad by\,\,putting\,x+t=x\\
&= \int_{\mathbb{R}}e^{-\frac{i}{2}(\frac{a}{b})t^{2}}\widehat{g}(x)\widehat{f}(x+t) dx \qquad when\,\,t\,\,is\,replace\,by\,-t\\
&=(g\otimes f)(t).
\end{align*}
(ii) Associativity
\begin{align*}
&[f\otimes (g\otimes h)](t)= \int_{\mathbb{R}}e^{-\frac{i}{2}(\frac{a}{b})t^{2}}\widehat{f}(x)\widehat{(g\otimes h)}(x+t) dx\\
&= \int_{\mathbb{R}}e^{-\frac{i}{2}(\frac{a}{b})t^{2}}\widehat{f}(x)e^{\frac{i}{2}(\frac{a}{b})(x+t)^{2}}(h\otimes g)(x+t) dx\\
&= \int_{\mathbb{R}}e^{-\frac{i}{2}(\frac{a}{b})t^{2}}\widehat{f}(x)e^{\frac{i}{2}(\frac{a}{b})(x+t)^{2}}\int_{\mathbb{R}}e^{-\frac{i}{2}(\frac{a}{b})(x+t)^{2}}\widehat{h}(y)\widehat{g}(y+x+t) dy  dx\\
&= \int_{\mathbb{R}}e^{-\frac{i}{2}(\frac{a}{b})t^{2}}\widehat{h}(y)\int_{\mathbb{R}}\widehat{f}(x)\widehat{g}(x+y+t) dx dy,\quad by\, changing\, order,\\
&= \int_{\mathbb{R}}e^{-\frac{i}{2}(\frac{a}{b})t^{2}}\widehat{h}(y)e^{\frac{i}{2}(\frac{a}{b})(y+t)^{2}}\int_{\mathbb{R}}e^{-\frac{i}{2}(\frac{a}{b})(y+t)^{2}}\widehat{f}(x)\widehat{g}(x+y+t) dx  dy\\
&= \int_{\mathbb{R}}e^{-\frac{i}{2}(\frac{a}{b})t^{2}}\widehat{h}(y)e^{\frac{i}{2}(\frac{a}{b})(y+t)^{2}}(f\otimes g)(y+t)dy\\
&= \int_{\mathbb{R}}e^{-\frac{i}{2}(\frac{a}{b})t^{2}}\widehat{h}(y)\widehat{(f\otimes g)}(y+t)dy=[(f\otimes g)\otimes h](t).
\end{align*}
\begin{lemma}
The space $(\mathcal{L}^{1}(\mathbb{R}), \Theta) $  is commutative semigroup.
\end{lemma}
\textbf{Proof:} Let $f, g, h\in \mathcal{L}^{1}(\mathbb{R})$, we have
(i)  commutative property
\begin{align*}
&(f\Theta g)(t)= \int_{-\infty}^{\infty}\frac{e^{-\frac{i}{2}(\frac{a}{b})t^{2}}}{2}\widehat{f}(x)[\widehat{g}(|x-t|)-\widehat{g}(x+t)] dx\\
&\leq \frac{1}{2}\int_{-\infty}^{\infty}e^{-\frac{i}{2}(\frac{a}{b})t^{2}}e^{\frac{i}{2}(\frac{a}{b})x^{2}}f(x)e^{\frac{i}{2}(\frac{a}{b})(x-t)^{2}}g(|x-t|)dx\\
&\quad- \frac{1}{2}\int_{-\infty}^{\infty}e^{-\frac{i}{2}(\frac{a}{b})t^{2}}e^{\frac{i}{2}(\frac{a}{b})x^{2}}f(x)e^{\frac{i}{2}(\frac{a}{b})(x+t)^{2}}g(x+t) dx = I_{1} + I_{2},
\end{align*}
we put $x-t=y$ in $I_{1}$ and $x+t=y$ in $I_{2}$, we get
\begin{align*}
&\leq \frac{1}{2}\int_{-\infty}^{\infty}e^{-\frac{i}{2}(\frac{a}{b})t^{2}}e^{\frac{i}{2}(\frac{a}{b})(y+t)^{2}}f(y+t)e^{\frac{i}{2}(\frac{a}{b})y^{2}}g(y)dy\\
&\quad- \frac{1}{2}\int_{-\infty}^{\infty}e^{-\frac{i}{2}(\frac{a}{b})t^{2}}e^{\frac{i}{2}(\frac{a}{b})(y-t)^{2}}f(y-t)e^{\frac{i}{2}(\frac{a}{b})y^{2}}g(y) dy \\
&=\int_{-\infty}^{\infty}\frac{e^{-\frac{i}{2}(\frac{a}{b})t^{2}}}{2}[\widehat{f}(y+t)-\widehat{f}(|y-t|)]\widehat{g}(y)]dy\\
&=\int_{\mathbb{R}}\frac{e^{-\frac{i}{2}(\frac{a}{b})t^{2}}}{2}\widehat{g}(y)[\widehat{f}(|y-t|)-\widehat{f}(y+t)]dy = (g\Theta f)(t).
\end{align*}
(ii) Associativity
 \begin{align*}
&[f\Theta (g\Theta h)](x)= \int_{-\infty}^{\infty}\frac{e^{-\frac{i}{2}(\frac{a}{b})x^{2}}}{2}[\widehat{f}(|y-x|)-\widehat{f}(y+x)] \widehat{(g\star h)}(y)dy\\
&= \int_{-\infty}^{\infty}\frac{e^{-\frac{i}{2}(\frac{a}{b})x^{2}}}{2}e^{\frac{i}{2}(\frac{a}{b})y^{2}}[\widehat{f}(y+x)+\widehat{f}(|y-x|)]  (g\star h)(y)dy\\
&= \int_{-\infty}^{\infty}\frac{e^{-\frac{i}{2}(\frac{a}{b})x^{2}}}{2}e^{\frac{i}{2}(\frac{a}{b})y^{2}}[\widehat{f}(|y-x|)-\widehat{f}(y+x)] \\
&\quad\times\bigg{(}\int_{-\infty}^{\infty}\frac{e^{-\frac{i}{2}(\frac{a}{b})y^{2}}}{2}\widehat{h}(z)[\widehat{g}(|z-y|)-\widehat{g}(z+y)] dz \bigg{)}dy\\
&= \int_{-\infty}^{\infty}\frac{e^{-\frac{i}{2}(\frac{a}{b})x^{2}}}{2}\widehat{h}(z) \bigg{(}\int_{-\infty}^{\infty}\frac{e^{-\frac{i}{2}(\frac{a}{b})y^{2}}}{2}e^{\frac{i}{2}(\frac{a}{b})y^{2}}\\
&\quad\times[\widehat{f}(|y-x|)-\widehat{f}(y+x)][\widehat{g}(|z-y|)-\widehat{g}(z+y)] dz \bigg{)}dy\\
&= \int_{-\infty}^{\infty}\frac{e^{-\frac{i}{2}(\frac{a}{b})x^{2}}}{2}\widehat{h}(z)\bigg{(}\int_{-\infty}^{\infty} \frac{e^{-\frac{i}{2}(\frac{a}{b})y^{2}}}{2}e^{\frac{i}{2}(\frac{a}{b})y^{2}}
[\widehat{f}(|y-x|)-\widehat{f}(y+x)]\widehat{g}(|z-y|)dy\\
&\quad+\int_{-\infty}^{\infty} \frac{e^{-\frac{i}{2}(\frac{a}{b})y^{2}}}{2}e^{\frac{i}{2}(\frac{a}{b})y^{2}}[\widehat{f}(|y-x|)-\widehat{f}(y+x)]\widehat{g}(z+y)dy \bigg{)} dz.
\end{align*}
Here we put $y=(u-z)$ in first integration and $y=(u+z)$ in second integration, we get,
 \begin{align*}
&= \int_{-\infty}^{\infty}\frac{e^{-\frac{i}{2}(\frac{a}{b})x^{2}}}{2}\widehat{h}(z)\\
&\quad\bigg{(}\int_{-\infty}^{\infty} \frac{e^{-\frac{i}{2}(\frac{a}{b})(u-z)^{2}}}{2}e^{\frac{i}{2}(\frac{a}{b})(u-z)^{2}}[\widehat{f}(|u-z-x|)-\widehat{f}(u-z+x)]\widehat{g}(|u|) du\\
&\quad-\int_{-\infty}^{\infty} \frac{e^{-\frac{i}{2}(\frac{a}{b})(u+z)^{2}}}{2}e^{\frac{i}{2}(\frac{a}{b})(u+z)^{2}}[\widehat{f}(|u+z-x|)-\widehat{f}(u+z+x)]\widehat{g}(u)du \bigg{)} dz\\
&= \int_{-\infty}^{\infty}\frac{e^{-\frac{i}{2}(\frac{a}{b})x^{2}}}{2}\widehat{h}(z)[e^{\frac{i}{2}(\frac{a}{b})(u-z)^{2}}(f\Theta g)(u-z)+e^{\frac{i}{2}(\frac{a}{b})(u+z)^{2}}(f\Theta g)(u+z)]dz\\
&= \int_{-\infty}^{\infty}\frac{e^{-\frac{i}{2}(\frac{a}{b})x^{2}}}{2}\widehat{h}(z)[\widehat{(f\Theta g)}(u-z)+\widehat{(f\Theta g)}(u+z)]dz.
\end{align*}
 consider $(u-z)=(x-z)$ in first integration and $(u+z)=(x+z)$ in second integration, we get,
 \begin{align*}
&= \int_{-\infty}^{\infty}\frac{e^{-\frac{i}{2}(\frac{a}{b})x^{2}}}{2}\widehat{h}(z)[\widehat{(f\Theta g)}(x-z)+\widehat{(f\Theta g)}(x+z)]dz=[(f\Theta g)\Theta h](x).
\end{align*}
 OR \\
 Now we prove this theorem by another way. We have, the operation $\star$ and $\otimes$ are commutative and associative, therefore\\
(i) commutativity
\begin{align*}
(f\Theta g)(t)= (f\star g)-(f\otimes g)=(g\star f)-(g\otimes f)=(g\Theta f)(t)
\end{align*}
(ii) Associativity
\begin{align*}
[(f\Theta g)\Theta h](t)&=[(f\star g)\star h](t)-[(f\otimes g)\otimes h](t)\\
&=[f\star (g\star h)](t)-[f\otimes (g\otimes h)](t)
\end{align*}
 Hence the space $(\mathcal{L}^{1}(\mathbb{R}), \Theta)$  is commutative semigroup.\hfill $\blacksquare$
\begin{theorem}\label{thm3.4.1}
 Let $f, g\in \mathcal{L}^{1} (\mathbb{R})$ and $F_{s}(s), G_{c}(s)$ denote the canonical sine and canonical cosine transform of $f(t)$ and $g(t)$ respectively, then
\begin{align}
\mathcal{S}_{A}[(f\Theta g)(t)](s)=\sqrt{2\pi ib}\,e^{-\frac{i}{2}(\frac{d}{b})s^{2}} F_{s}(s)G_{c}(s).
\end{align}
\end{theorem}
\textbf{Proof:} Let $f, g\in \mathcal{L}^{1} (\mathbb{R})$ then we have,
\begin{align*}
&\mathcal{S}_{A}[(f\Theta g)(t)](s)=\frac{e^{\frac{i}{2}(\frac{d}{b})s^{2}}}{\sqrt{2\pi ib}}\int_{-\infty}^{\infty}sin(\frac{s}{b}t)e^{\frac{i}{2}(\frac{a}{b})t^{2}}\\
&\quad\times\bigg{(}\int_{-\infty}^{\infty}\frac{e^{-\frac{i}{2}(\frac{a}{b})t^{2}}}{2}\widehat{f}(x)[\widehat{g}(|x-t|)-\widehat{g}(x+t)] dx \bigg{)}dt\\
&=\frac{e^{\frac{i}{2}(\frac{d}{b})s^{2}}}{2\sqrt{2\pi ib}}\int_{-\infty}^{\infty}\int_{-\infty}^{\infty}e^{\frac{i}{2}(\frac{a}{b})x^{2}}f(x)\\
&\quad\times[e^{\frac{i}{2}(\frac{a}{b})(x-t)^{2}}g(x-t)- e^{\frac{i}{2}(\frac{a}{b})(x+t)^{2}}g(x+t)]sin(\frac{s}{b}t)dt dx\\
&=\frac{e^{\frac{i}{2}(\frac{d}{b})s^{2}}}{2\sqrt{2\pi ib}}\int_{-\infty}^{\infty}\int_{-\infty}^{\infty}e^{\frac{i}{2}(\frac{a}{b})x^{2}}f(x)e^{\frac{i}{2}(\frac{a}{b})(x-t)^{2}}g(x-t)sin(\frac{s}{b}t)dt dx\\
&\quad -\frac{e^{\frac{i}{2}(\frac{d}{b})s^{2}}}{2\sqrt{2\pi ib}}\int_{-\infty}^{\infty}\int_{-\infty}^{\infty}e^{\frac{i}{2}(\frac{a}{b})x^{2}}f(x)e^{\frac{i}{2}(\frac{a}{b})(x+t)^{2}}g(x+t)sin(\frac{s}{b}t)dt dx\\
&= I_{1}-I_{2},
\end{align*}
we put $x-t=y$ in $I_{1}$ and $x+t=y$ in $I_{2}$, we get
\begin{align*}
&=\frac{e^{\frac{i}{2}(\frac{d}{b})s^{2}}}{2\sqrt{2\pi ib}}\int_{-\infty}^{\infty}\int_{-\infty}^{\infty}e^{\frac{i}{2}(\frac{a}{b})x^{2}}f(x)e^{\frac{i}{2}(\frac{a}{b})(y)^{2}}g(y)sin(\frac{s}{b}(y+x))dy dx\\
&\quad -\frac{e^{\frac{i}{2}(\frac{d}{b})s^{2}}}{2\sqrt{2\pi ib}}\int_{-\infty}^{\infty}\int_{-\infty}^{\infty}e^{\frac{i}{2}(\frac{a}{b})x^{2}}f(x)e^{\frac{i}{2}(\frac{a}{b})(y)^{2}}g(y)sin(\frac{s}{b}(y-x))dy dx\\
&=\frac{e^{\frac{i}{2}(\frac{d}{b})s^{2}}}{2\sqrt{2\pi ib}}\int_{-\infty}^{\infty}\int_{-\infty}^{\infty}e^{\frac{i}{2}(\frac{a}{b})x^{2}}e^{\frac{i}{2}(\frac{a}{b})(y)^{2}}\\
&\quad\times [sin(\frac{s}{b}(y+x))-sin(\frac{s}{b}(y-x))]g(y)f(x)dy dx\\
&=\frac{e^{\frac{i}{2}(\frac{d}{b})s^{2}}}{2\sqrt{2\pi ib}}\\
&\quad\times\int_{-\infty}^{\infty}\int_{-\infty}^{\infty}e^{\frac{i}{2}(\frac{a}{b})x^{2}}e^{\frac{i}{2}(\frac{a}{b})(y)^{2}}[2cos(\frac{s}{b}y)sin(\frac{s}{b}x)]g(y)f(x)dy dx\\
&=\frac{\sqrt{2\pi ib}}{e^{\frac{i}{2}(\frac{d}{b})s^{2}}}\frac{e^{\frac{i}{2}(\frac{d}{b})s^{2}}}{\sqrt{2\pi ib}}\int_{-\infty}^{\infty}e^{\frac{i}{2}(\frac{a}{b})x^{2}}sin(\frac{s}{b}x)f(x)G_{c}(s)dx\\
&=\sqrt{2\pi ib}\,e^{-\frac{i}{2}(\frac{d}{b})s^{2}}F_{s}(s)G_{c}(s).
\end{align*}$\hfill\blacksquare$

\begin{lemma}\label{lemma3.4.1} Let $f, g\in \mathcal{L}^{1}(\mathbb{R})$, then $\|(f\Theta g)\|\leq\|f\|_{1} \|g\|_{1}$ and hence the convolution is a member of $\mathcal{L}^{1} (\mathbb{R})$.
\end{lemma}
\textbf{Proof:} Let $f, g\in \mathcal{L}^{1}(\mathbb{R})$  we have,
\begin{align*}
&\|(f\Theta g)\|_{1}\\
&=\int_{\mathbb{R}}|(f\Theta g)|dt=\int_{\mathbb{R}}\bigg{|}\int_{-\infty}^{\infty}\frac{e^{-\frac{i}{2}(\frac{a}{b})t^{2}}}{2}\widehat{f}(x)[\widehat{g}(|x-t|)-\widehat{g}(x+t)] dx\bigg{|}dt\\
&\leq \frac{1}{2}\int_{\mathbb{R}}\bigg{|}\int_{-\infty}^{\infty}e^{-\frac{i}{2}(\frac{a}{b})t^{2}}e^{\frac{i}{2}(\frac{a}{b})x^{2}}f(x)e^{\frac{i}{2}(\frac{a}{b})(x-t)^{2}}g(|x-t|)dx\bigg{|}dt\\
&+\frac{1}{2}\int_{\mathbb{R}}\bigg{|}\int_{-\infty}^{\infty}e^{-\frac{i}{2}(\frac{a}{b})t^{2}}e^{\frac{i}{2}(\frac{a}{b})x^{2}}f(x)e^{\frac{i}{2}(\frac{a}{b})(x+t)^{2}}g(x+t) dx\bigg{|}dt\\
&\leq \frac{1}{2}\int_{\mathbb{R}}\int_{-\infty}^{\infty}|g(|x-t|)f(x)|dx dt+\frac{1}{2}\int_{\mathbb{R}}\int_{-\infty}^{\infty}|g(x+t)f(x)|dx dt\\
&\leq \frac{1}{2}\int_{\mathbb{R}}\|g(|x-t|)f(x)\|_{1}dt+\frac{1}{2}\int_{\mathbb{R}}\|g(x+t)f(x)\|_{1}dt\\
&\leq \frac{1}{2}\|f\|_{1}\int_{\mathbb{R}}|g(|x-t|)|dt+\frac{1}{2}\|f(x)\|_{1}\int_{\mathbb{R}}|g(x+t)|dt= \|f\|_{1} \|g\|_{1}.
\end{align*}$\hfill\blacksquare$
\begin{lemma}\label{lemma3.4.3}
  Let $( \varphi_{n} ), ( \psi_{n} )\in \bigtriangleup$ for $n\in \mathbb{N}$ then the convolution $(\varphi_{n}\Theta \psi_{n}) \in \bigtriangleup$ for each $n\in \mathbb{N}$ .
\end{lemma}
\textbf{Proof:} Let $( \varphi_{n} ), ( \psi_{n} )\in \bigtriangleup$ for $n\in \mathbb{N}$, we have\\
(i)
\begin{align*}
&\int_{\mathbb{R}}e^{(\frac{ia}{b})tx}(\varphi_{n}\Theta\psi_{n})(t)dt\\
&=\int_{\mathbb{R}}e^{(\frac{ia}{b})tx}\int_{-\infty}^{\infty}\frac{e^{-\frac{i}{2}(\frac{a}{b})t^{2}}}{2}\widehat{\varphi_{n}}(x)[\widehat{\psi_{n}}(|t-x|)-\widehat{\psi_{n}}(t+x)] dx dt\\
&=\frac{1}{2} \int_{\mathbb{R}}e^{(\frac{ia}{b})tx}\int_{-\infty}^{\infty}e^{-\frac{i}{2}(\frac{a}{b})t^{2}}e^{\frac{i}{2}(\frac{a}{b})x^{2}}\varphi_{n}(x)e^{\frac{i}{2}(\frac{a}{b})(t-x)^{2}}\psi_{n}(t-x) dxdt\\
&\quad- \frac{1}{2}\int_{\mathbb{R}}e^{(\frac{ia}{b})tx}\int_{-\infty}^{\infty}e^{-\frac{i}{2}(\frac{a}{b})t^{2}}e^{\frac{i}{2}(\frac{a}{b})x^{2}}\varphi_{n}(x)e^{\frac{i}{2}(\frac{a}{b})(t+x)^{2}}\psi_{n}(t+x)dx dt\\
&= \frac{1}{2}\int_{\mathbb{R}}e^{(\frac{ia}{b})tx}\int_{-\infty}^{\infty}e^{(\frac{ia}{b})x^{2}}\varphi_{n}(x)[e^{(\frac{-ia}{b})tx}\psi_{n}(t-x)-e^{(\frac{ia}{b})tx}\psi_{n}(t+x)] dxdt\\
&= \frac{1}{2}\int_{-\infty}^{\infty}e^{(\frac{ia}{b})tx}\varphi_{n}(x)e^{(\frac{ia}{b})x^{2}}\\
&\quad\times\bigg{[}\int_{\mathbb{R}}e^{(\frac{-ia}{b})tx}\psi_{n}(t-x) dt-\int_{\mathbb{R}}e^{(\frac{ia}{b})tx}\psi_{n}(t+x)dt\bigg{]} dx\\
&= \frac{1}{2}\int_{-\infty}^{\infty}e^{(\frac{ia}{b})tx}\varphi_{n}(x)e^{(\frac{ia}{b})x^{2}}\\
&\quad\times\bigg{[}\int_{\mathbb{R}}e^{(\frac{-ia}{b})x(y+x)}\psi_{n}(y) dy-\int_{\mathbb{R}}e^{(\frac{ia}{b})x(y-x)}\psi_{n}(y)dy\bigg{]} dx\\
&= \frac{1}{2}\int_{-\infty}^{\infty}e^{(\frac{ia}{b})tx}\varphi_{n}(x)\bigg{[}\int_{\mathbb{R}}e^{(\frac{-ia}{b})xy}\psi_{n}(y) dy-\int_{\mathbb{R}}e^{(\frac{ia}{b})xy}\psi_{n}(y)dy\bigg{]} dx\\
&= \frac{1}{2}\int_{-\infty}^{\infty}e^{(\frac{ia}{b})tx}\varphi_{n}(x)\bigg{[}\int_{\mathbb{R}}e^{(\frac{ia}{b})zy}\psi_{n}(y) dy+\int_{\infty}^{-\infty}e^{(\frac{ia}{b})xy}\psi_{n}(y)dy\bigg{]} dx=1.
\end{align*}
(ii) From lemma \eqref{lemma3.4.1}, we have
\begin{align*}
\|(\varphi_{n}\Theta\psi_{n})\|_{1}&=\int_{\mathbb{R}}|(\varphi_{n}\Theta\psi_{n})(t)|dt
\leq \|\varphi_{n}\|_{1}\|\psi_{n}\|_{1}\leq M \quad for \,\, M > 0 .
\end{align*}
(iii)
\begin{align*}
&\lim_{n\rightarrow \infty}\int_{|t|>\epsilon} |(\varphi_{n}\Theta\psi_{n})(t)| dt \leq \frac{1}{2}\lim_{n\rightarrow\infty}\int_{|t|>\epsilon}\\
&\quad\bigg{(} \int_{-\infty}^{\infty}\big{|}\varphi_{n}(x)\psi_{n}(|t-x|)\big{|}dx-\int_{-\infty}^{\infty}\big{|}\varphi_{n}(x)\psi_{n}(t+x)\big{|}dxdtdt\bigg{)}\\
&\leq \lim_{n\rightarrow \infty}\int_{|t|>\epsilon}\|\varphi_{n}(x)\psi_{n}(|t-x|)\|_{1}dt+\frac{1}{2}\lim_{n\rightarrow \infty}\int_{|t|>\epsilon}\|\varphi_{n}(x)\psi_{n}(t+x)\|_{1}dt\\
&\leq \frac{\|\varphi_{n}\|_{1}}{2}\lim_{n\rightarrow \infty}\int_{|t|>\epsilon}|\psi_{n}(|t-x|)\big{|}dt+\frac{\|\varphi_{n}\|_{1}}{2}\lim_{n\rightarrow \infty}\int_{|t|>\epsilon}|\psi_{n}(t+x)|dt\\
&\leq \|\varphi_{n}\|_{1}\lim_{n\rightarrow \infty}\int_{|y| > \delta}|\psi_{n}(y)\big{|}dy\longrightarrow 0 \quad for\, each\quad \delta > 0.
\end{align*}
\begin{lemma}\label{lemma3.4.4}
Let $f \in \mathcal{L}^{1}(\mathbb{R})$ and $( \delta_{n} ) \in \bigtriangleup$ then $f \Theta \delta_{n}\longrightarrow f$ as  $n\longrightarrow \infty $ in $\mathcal{L}^{1}(\mathbb{R})$.
\end{lemma}
\textbf{Proof:} Let $f\in \mathcal{L}^{1}(\mathbb{R})$ and $( \delta_{n} ) \in \bigtriangleup$, we have the relation $f \Theta \delta_{n}= f\star \delta_{n}- f \otimes \delta_{n}$ for each $n\in\mathbb{N}$. By lemma \eqref{lemma3.2.5} the first term of above equality tend to $f$ as $n\longrightarrow \infty$. Since $\int_{\mathbb{R}}|\delta_{n}(t)|dt\leq M$ for some $M>0$ and $supp\delta_{n}\longrightarrow\{0\}$ as $n\longrightarrow\infty$, implies $\int_{|t|>\epsilon}|\delta_{n}(t)|\rightarrow 0 \,\, as \,\,n\rightarrow\infty $, we have
\begin{align*}
\|f \otimes \delta_{n}\|_{1}&\leq\int_{\mathbb{R}}\bigg{|}\int_{\mathbb{R}}e^{-\frac{i}{2}(\frac{a}{b})t^{2}}\widehat{\delta_{n}}(x)\widehat{f}(x+y) dx\bigg{|}dt\\
&\leq \int_{\mathbb{R}}|\delta_{n}(t)|\int_{\mathbb{R}}|f(x+y)|dx dt\\
&=\|f\|_{1} \int_{|t|>\epsilon}|\delta_{n}(t)|\leq M_{1} \int_{|t|>\epsilon}|\delta_{n}(t)|\quad for \,\, M_{1}>0,
\end{align*}
Hence $f \Theta \delta_{n}\rightarrow f$ as  $n\rightarrow \infty $ in $\mathcal{L}^{1}(\mathbb{R})$.$\hfill\blacksquare$
\begin{lemma}\label{lemma3.4.6}
   Let $f_{n}\longrightarrow f$ as $\longrightarrow\infty$ in $\mathcal{L}^{1}(\mathbb{R})$ and $( \varphi_{n} )\in \bigtriangleup$  then $(f_{n}\Theta\varphi_{n})\rightarrow f$ as $n\rightarrow \infty$ in $\mathcal{L}^{1}(\mathbb{R})$.
\end{lemma}
\textbf{Proof:} Let $f_{n}\rightarrow f$ as $\rightarrow\infty$ in $\mathcal{L}^{1}(\mathbb{R})$ and $( \varphi_{n} )\in \bigtriangleup$, we have $\|f_{n}\Theta \varphi_{n}-f\|_{1}=\|f_{n}\Theta \varphi_{n}-f\Theta \varphi_{n}+f\Theta \varphi_{n}-f\|_{1}\leq \|f_{n}\Theta \varphi_{n}-f\Theta \varphi_{n}\|_{1}+\|f\Theta \varphi_{n}-f\|_{1}$. From lemma \eqref{lemma3.4.4} $\|f\Theta \varphi_{n}-f\|_{1}\rightarrow 0$, then we have $\|f_{n}\Theta \varphi_{n}-f\|_{1}\leq \|f_{n}\Theta \varphi_{n}-f\Theta \varphi_{n}\|_{1}\leq \|f_{n}-f\|_{1}\|\varphi_{n}\|_{1}\leq M\|f_{n}-f\|_{1}\rightarrow 0$ as $n\rightarrow \infty$. where $\|\varphi_{n}\|_{1}\leq M$ for $M > 0$.$\hfill\blacksquare$
\begin{lemma}\label{lemma3.4.8}
Let $( \delta_{n} )\in \bigtriangleup$ for all $n\in\mathbb{N}$ then $\mathcal{S}_{A}(\delta_{n})$ converges uniformly on each compact set to a constant function $1$ in $\mathcal{L}^{1}(\mathbb{R})$.
\end{lemma}
\textbf{Proof:} Let $( \delta_{n} ) \in \bigtriangleup$. Using property (i), (ii) of $\bigtriangleup$ we have  $supp(\delta_{n})\rightarrow 0$ for $n\rightarrow \infty$ on each compact subset $K$ of $\mathbb{R}$ and there exist $M > 0$, $\bigg{|}sin(\frac{s}{b}t)e^{\frac{i}{2}(\frac{a}{b})t^{2}}-\sqrt{2\pi ib}\,e^{\frac{-i}{2}(\frac{d}{b})s^{2}}e^{(\frac{ia}{b})xt}\bigg{|} < M $ for all $ t\in\mathbb{R}$, therefore
\begin{align*}
&\|\big{(}\mathcal{S}_{A}(\delta_{n})-1\big{)}\|_{1}\\
&=\int_{\mathbb{R}}\bigg{|} \frac{e^{\frac{i}{2}(\frac{d}{b})s^{2}}}{\sqrt{2\pi ib}}\int_{-\infty}^{\infty}sin(\frac{s}{b}t)e^{\frac{i}{2}(\frac{a}{b})t^{2}}\delta_{n}(t)dt-\int_{-\infty}^{\infty}e^{(\frac{ia}{b})xt}\delta_{n}(t)dt\bigg{|}ds\\
&=\int_{\mathbb{R}}\bigg{|} \frac{e^{\frac{i}{2}(\frac{d}{b})s^{2}}}{\sqrt{2\pi ib}}\int_{-\infty}^{\infty}\bigg{(}sin(\frac{s}{b}t)e^{\frac{i}{2}(\frac{a}{b})t^{2}}-\sqrt{2\pi ib}\,e^{\frac{-i}{2}(\frac{d}{b})s^{2}}e^{(\frac{ia}{b})xt}\bigg{)}\delta_{n}(t)dt\bigg{|}ds\\
&\leq \frac{M^{2}}{2\pi b}\int_{\mathbb{R}}\bigg{(}\int_{K}|\delta_{n}(t)|dt\bigg{)} ds \rightarrow 0 \quad as\,\, n\rightarrow \infty.
\end{align*} $\hfill\blacksquare$

\section{CST For Boehmians}

\begin{definition}  A pair of sequences $(( f_{n} ), ( \varphi_{n} ))$ is called a {\it quotient of sequences} , denoted by $f_{n}/\varphi_{n}$, if $( f_{n} )\in \mathcal{L}^{1}(\mathbb{R})$ , $( \varphi_{n} )\in\bigtriangleup$ and $f_{m}\Theta\varphi_{n}=f_{n}\Theta\varphi_{m}$ $\forall\, m, n\in \mathbb{N}$.
\end{definition}
\begin{definition}
Let $( f_{n} ), \{g_{n}\}\in \mathcal{L}^{1}(\mathbb{R})$ and $( \psi_{n} ), ( \varphi_{n} )\in\bigtriangleup$. Two quotient of sequences $f_{n}/ \varphi_{n}$ and $g_{n} / \psi_{n}$ are equivalent if $f_{n}\Theta\psi_{n}=g_{n}\Theta\varphi_{n}$ $\forall\, n\in\mathbb{N}$. This is an equivalence relation.
\end{definition}
\begin{definition}
 The equivalence class of quotient of sequences is called a {\it Boehmian} and space of all Boehmians denoted by $\mathcal{B}_{\Theta}=\mathcal{B'}(\mathcal{L}^{1}(\mathbb{R}), \bigtriangleup, \Theta)$. Let function $f\in \mathcal{L}^{1}(\mathbb{R})$ then it can be identified with the Boehmian $[f\Theta\delta_{n}/\delta_{n}]$, where $( \delta_{n} )\in\bigtriangleup$. Let the Boehmian $F=[f_{n}/\delta_{n}]$, then $F\Theta\delta_{n}=f_{n} \in \mathcal{L}^{1}(\mathbb{R})$ $\forall\, n\in\mathbb{N}$.
\end{definition}
\begin{definition}
A sequence of Boehmians $( f_{n} )$ is called $\Delta-$convergent to a Boehmian $F$ ($\Delta-\lim F_{n}=F$) if there exist a delta sequence $( \delta_{n} )$ such that  $(F_{n}-F)\Theta\delta_{n}\in \mathcal{L}^{1}(\mathbb{R})$, for every $n\in\mathbb{N}$ and that $\|((F_{n}-F)\Theta\delta_{n})\|_{1}\rightarrow 0$ as $n\rightarrow \infty$.
\end{definition}
\begin{definition}
 A sequence of Boehmians $( f_{n} )$ is called $\delta-$ convergent to a Boehmian $F$ ($\delta-\lim F_{n}=F$) if there exist a delta sequence $( \delta_{n} )$ such that  $F_{n}\Theta\delta_{k}\in \mathcal{L}^{1}(\mathbb{R})$ and $F\Theta\delta_{k}\in \mathcal{L}^{1}(\mathbb{R})$ for every $n, k\in\mathbb{N}$ and that $\|((F_{n}-F)\Theta\delta_{k})\|_{1}\rightarrow 0$ as $n\rightarrow \infty$ for each $k\in\mathbb{N}$.
 \end{definition}
\indent\quad Let $( \delta_{n} )$ is a delta sequence, then $[\delta_{n}/\delta_{n}]$ represents a Boehmian for each $n\in\mathbb{N}$. Since the Boehmian [$\delta_{n}/\delta_{n}$] correspond to Dirac delta distribution $\delta$, all the derivative of $\delta$ are also Boehmian. Further, if $(\delta_{n})$ is infinitely differentiable and bounded support, then the $k^{th}$ derivative of $\delta$ is define as $\delta^{(k)}=[\delta_{n}^{(k)}/\delta_{n}]$, implies  $\delta^{(k)}\in \mathcal{B}_{\Theta}$, where $k\in\mathbb{N}$. The $k^{th}$ derivative of Boehmian $F\in \mathcal{B}_{\Theta}$ is define by $F^{(k)}=F\Theta\delta^{(k)}$.\\
The scalar multiplication, addition and convolution in $\mathcal{B}_{\Theta}$ is define as below,
\begin{align*}
\lambda[f_{n}/\varphi_{n}]&=[\lambda f_{n}/\varphi_{n}]\\
[f_{n}/\varphi_{n}]+[g_{n}/\psi_{n}]&=[(f_{n}\Theta\psi_{n}+g_{n}\Theta\varphi_{n})/\varphi_{n}\Theta\psi_{n}]\\
[f_{n}/\varphi_{n}]\Theta[g_{n}/\psi_{n}]&=[f_{n}\Theta g_{n}/\varphi_{n}\Theta\psi_{n}]
\end{align*}
\begin{definition}
Let $( f_{n} )\in \mathcal{L}^{1}(\mathbb{R})$ and $( \delta_{n} )\in\bigtriangleup$ for $n\in\mathbb{N}$, we define the canonical sine transform $\mathcal{S}_{A} : \mathcal{B}_{\Theta}\rightarrow\mathcal{B}$ as
 \begin{align}
 \mathcal{S}_{A} [f_{n}/\delta_{n}]= \mathcal{S}_{A} (f_{n})/ \mathcal{C}_{A} (\delta_{n}) \qquad for \quad [f_{n}/\delta_{n}]\in \mathcal{B}_{\Theta}.
\end{align}
\end{definition}
\indent\quad The CST on $\mathcal{B}_{\Theta}$ is well defined. Indeed if  $[f_{n}/\delta_{n}]\in \mathcal{B}_{\Theta}$, then $f_{n}\Theta\delta_{m}=f_{m}\Theta\delta_{n}$ for all $m, n \in \mathbb{N}$. Applying the CST on both sides, we get $\mathcal{S}_{A} (f_{n}) \mathcal{C}_{A} (\delta_{m})=\mathcal{S}_{A} (f_{m}) \mathcal{C}_{A} (\delta_{n})$ for all $m, n \in \mathbb{N}$ and hence $ \mathcal{S}_{A} (f_{n})/ \mathcal{C}_{A} (\delta_{n}) \in \mathcal{B} $. Further if $[f_{n}/\psi_{n}]=[g_{n}/\delta_{n}]\in \mathcal{B}_{\Theta}$ then we have $f_{n}\Theta\delta_{n}=g_{n}\Theta\psi_{n}$ for all $ n \in \mathbb{N}$. Again applying the CST on both sides, we get $\mathcal{S}_{A} (f_{n}) \mathcal{C}_{A} (\delta_{n})=\mathcal{S}_{A} (g_{n}) \mathcal{C}_{A} (\psi_{n})$ for all $ n \in \mathbb{N}$. i.e. $\mathcal{S}_{A} (f_{n})/ \mathcal{C}_{A} (\psi_{n})=\mathcal{S}_{A} (g_{n})/ \mathcal{C}_{A} (\delta_{n})$ in $\mathcal{B}$.

\begin{lemma}
Let $( f_{n} )\in \mathcal{L}^{1}(\mathbb{R})$ for $n\in\mathbb{N}$, then
\begin{align}
\mathcal{S}_{A} [f_{n}](s)= \sqrt{\frac{1}{2\pi ib}}e^{\frac{i}{2}(\frac{d}{b})s^{2}}\int_{-\infty}^{\infty}cos(\frac{s}{b}t)e^{\frac{i}{2}(\frac{a}{b})t^{2}}f_{n}(t)dt
\end{align}
converges uniformly on each compact subset of $\mathbb{R}$.
\end{lemma}
\textbf{Proof:} Let $( \delta_{n} )$ is delta sequence, then $\mathcal{C}_{A}[\delta_{n}]$ converges uniformly on each compact set to a function $1$. Hence, for each compact subset $M$, $\mathcal{S}_{A}[\delta_{m}] > 0$  on $M$, for almost all $m\in\mathbb{N}$ Moreover,
\begin{align}
\mathcal{S}_{A} (f_{n})&= \mathcal{S}_{A} (f_{n})\frac{\mathcal{C}_{A} (\delta_{m})}{\mathcal{C}_{A} (\delta_{m})}=\frac{e^{\frac{i}{2}(\frac{d}{b})s^{2}}}{\sqrt{2\pi ib}}\frac{\mathcal{S}_{A}(f_{n}\Theta\delta_{m})}{\mathcal{C}_{A} (\delta_{m})}\nonumber\\
&=\frac{e^{\frac{i}{2}(\frac{d}{b})s^{2}}}{\sqrt{2\pi ib}}\frac{\mathcal{S}_{A}(f_{m}\Theta\delta_{n})}{\mathcal{C}_{A} (\delta_{m})}=\frac{\mathcal{S}_{A}(f_{m})}{\mathcal{C}_{A}(\delta_{m})}{\mathcal{C}_{A} (\delta_{n})},\qquad on\quad M,
\end{align}
as $n\rightarrow \infty$ we get $\mathcal{S}_{A} (f_{n})\rightarrow \frac{\mathcal{S}_{A}(f_{m})}{\mathcal{C}_{A}(\delta_{m})} $, on each compact subset $M$ of $\mathbb{R}$.$\hfill\blacksquare$

In view of the above lemma, the canonical sine transform of Boehmian $F\in \mathcal{B}_{\Theta}$  can be define as\\
\begin{align}
F=\lim_{n\rightarrow\infty}\mathcal{S}_{A} (f_{n}).
\end{align}

Now we have to show that the above definition is well defined. Let two Boehmians $F=[f_{n}/ \varphi_{n}]$ and $G=[g_{n} / \psi_{n}]$ are the equivalent. i.e. $f_{n} \Theta \psi_{n}=g_{n} \Theta \varphi_{n}$ for all $n\in\mathbb{N}$. Taking canonical sine transform on both side, we have
\begin{align*}
 \mathcal{S}_{A}[f_{n}\Theta\psi_{n}](s)&=\mathcal{S}_{A}[g_{n}\Theta\varphi_{n}](s)\\
\mathcal{S}_{A} [f_{n}](s)\mathcal{C}_{A}[\psi_{n}](s)&= \mathcal{S}_{A} [g_{n}](s)\mathcal{C}_{A}[\varphi_{n}](s)
\end{align*}
which implies that
\begin{align*}
\lim_{n\rightarrow \infty}\ \mathcal{S}_{A}[f_{n}](s)&=\lim_{n\rightarrow \infty}\mathcal{S}_{A}[g_{n}](s)\\
\mathcal{S}_{A}[F](s)&=\mathcal{S}_{A}[G](s).
\end{align*}
\begin{theorem}
The canonical sine transform $\mathcal{S}_{A}:\mathcal{B}_{\Theta}\longrightarrow \mathcal{B}$ is consistent with $\mathcal{S}_{A}: \mathcal{L}^{1}(\mathbb{R})\longrightarrow \mathcal{L}^{1}(\mathbb{R})$.
\end{theorem}
\textbf{Proof:} Let $f\in \mathcal{L}^{1}(\mathbb{R})$. The Boehmian representing $f$ in $\mathcal{B}_{\Theta}$ is $[(f\Theta\delta_{n})/\delta_{n}]$ where $( \delta_{n} )\in \bigtriangleup$. Then for each $n\in\mathbb{N}$, $\mathcal{S}_{A}((f\Theta\delta_{n})/\delta_{n})= \sqrt{2\pi ib}\,e^{-\frac{i}{2}(\frac{d}{b})s^{2}} \mathcal{S}_{A} (f) \mathcal{C}_{A} (\delta_{n})/\mathcal{C}_{A} (\delta_{n})$, which is a Boehmian in $\mathcal{B}$, which represents $\sqrt{2\pi ib} \,e^{-\frac{i}{2}(\frac{d}{b})s^{2}}\mathcal{S}_{A}(f)$.$\hfill\blacksquare$
\begin{theorem}
The canonical sine transform $\mathcal{S}_{A}:\mathcal{B}_{\Theta}\longrightarrow \mathcal{B}$ is a bijection.
\end{theorem}
\textbf{Proof:} Let $F=[f_{n}/\varphi_{n}], G=[g_{m}/\psi_{m}]\in \mathcal{B}_{\Theta}$ such that $\mathcal{S}_{A}(F)=\mathcal{S}_{A}(G)$. From this we get $\mathcal{S}_{A}(f_{n}) \mathcal{C}_{A}(\psi_{m})=\mathcal{S}_{A}(g_{m}) \mathcal{C}_{A}(\varphi_{n})$ for all $m, n\in\mathbb{N}$ and hence $\mathcal{S}_{A}(f_{n}\Theta\psi_{m})=\mathcal{S}_{A}(g_{m}\Theta\varphi_{n})$ for all $m, n\in\mathbb{N}$. Since CST is one-to-one in $\mathcal{L}^{1}(\mathbb{R})$, we get
$(f_{n}\Theta\psi_{m})=(g_{m}\Theta\varphi_{n})$ for all $m, n\in\mathbb{N}$. This implies $F=G$.
\par Let $G=[g_{m}/\psi_{m}]\in \mathcal{B}_{\Theta}$. Since $\mathcal{S}_{A}: \mathcal{L}^{1}(\mathbb{R})\longrightarrow \mathcal{L}^{1}(\mathbb{R})$ is onto. Choose $f_{n}\in \mathcal{L}^{1}(\mathbb{R})$ such that $g_{n}= \mathcal{S}_{A} (f_{n})$. Now using the relation $g_{n}\cdot \mathcal{C}_{A}(\psi_{m})=g_{m}\cdot \mathcal{C}_{A}(\psi_{n})$ for all $m, n\in\mathbb{N}$. we obtain $\mathcal{S}_{A}(f_{n})\cdot \mathcal{C}_{A}(\psi_{m})=\mathcal{S}_{A} (f_{m})\cdot \mathcal{C}_{A}(\psi_{n}) \Rightarrow \mathcal{S}_{A}(f_{n}\Theta\psi_{m})=\mathcal{S}_{A} (f_{m}\Theta\psi_{n})$ Since CST is one-to-one in $\mathcal{L}^{1}(\mathbb{R})$, we get $(f_{n}\Theta\psi_{m})=(f_{m}\Theta\psi_{n})$. Now if we take $F=[f_{n}/\psi_{n}]$ then $F\in \mathcal{B}_{\Theta}$ and $\mathcal{S}_{A}(F)=G$. Thus the theorem is hold.$\hfill\blacksquare$

\begin{theorem}
Let $F, G \in \mathcal{B}_{\Theta}$ then\\
(1) $\mathcal{S}_{A}(\lambda F + G)= \lambda \mathcal{S}_{A} (F)+ \mathcal{S}_{A} (G)$, for any complex $\lambda$,\\
(2) $\mathcal{S}_{A} (F\Theta G)= \mathcal{S}_{A}(F) \mathcal{S}_{A}(G)$,\\
(3) $\mathcal{S}_{A} (F(kt))=\frac{1}{k}e^{(1-\frac{1}{k^{2}})\frac{i}{2}\frac{d}{b}s^{2}}\mathcal{S}_{A} (e^{(\frac{1}{k^{2}}-1)\frac{i}{2}\frac{a}{b}t^{2}}F)(\frac{s}{k})$\\
(4) $\mathcal{S}_{A}(F(x+t))=e^{\frac{i}{2}(\frac{a}{b})\tau^{2}}\{cos(\frac{s}{b}\tau)\mathcal{S}_{A}(e^{-it\tau(\frac{a}{b})} F)-sin(\frac{s}{b}\tau)\mathcal{C}_{A}(e^{-it\tau(\frac{a}{b})} F)\}$\\
(5) $\mathcal{S}_{A}(e^{ixt} F)(s)=\mathcal{S}_{A}(cosxt F)(s)+i \mathcal{S}_{A} (sinxt F)(s)$\\
(6) $\mathcal{S}_{A}\big{(}cos xt\, F\big{)}(s)=\frac{e^{\frac{i}{2}(db)x^{2}}}{2}\big{(}e^{idxs}\mathcal{S}_{A} \big{(}F\big{)}(s+bx)+e^{-idxs}\mathcal{S}_{A} \big{(}F\big{)}(s-bx)\big{)} $\\
(7) $\mathcal{S}_{A}\big{(}sin xt\, F\big{)}(s)=\frac{e^{\frac{i}{2}(db)x^{2}}}{2}\big{(}e^{-idxs}\mathcal{C}_{A} \big{(}F\big{)}(s-bx)-e^{idxs}\mathcal{C}_{A} \big{(}F\big{)}(s+bx)\big{)} $.
\end{theorem}
\textbf{Proof:} Let $F, G \in \mathcal{B}_{\Theta}$ and $\lambda $ be any scalar,\\
(1)Linearity
\begin{align*}
\mathcal{S}_{A}[(\lambda F + G)(t)](s)&=\frac{1}{\sqrt{2\pi ib}}e^{\frac{i}{2}(\frac{d}{b})s^{2}}\int_{-\infty}^{\infty}sin(\frac{s}{b}t)e^{\frac{i}{2}(\frac{a}{b})t^{2}}(\lambda F + G)(t)dt\\
&=\frac{1}{\sqrt{2\pi ib}}e^{\frac{i}{2}(\frac{d}{b})s^{2}}\lambda\int_{-\infty}^{\infty}sin(\frac{s}{b}t)e^{\frac{i}{2}(\frac{a}{b})t^{2}} F(t)dt\\
&\quad+\frac{1}{\sqrt{2\pi ib}}e^{\frac{i}{2}(\frac{d}{b})s^{2}}\int_{-\infty}^{\infty}sin(\frac{s}{b}t)e^{\frac{i}{2}(\frac{a}{b})t^{2}}G(t)dt\\
&=\lambda\mathcal{S}_{A}[ F (t)](s)+\mathcal{S}_{A}[G(t)](s).
\end{align*}
(2) Convolution of Boehmians
\begin{align*}
\mathcal{S}_{A} (F \Theta G)&=\mathcal{S}_{A} ([f_{n}/\varphi_{n}]\Theta[g_{n}/\psi_{n}])=\mathcal{S}_{A} [(f_{n}\Theta g_{n})/(\varphi_{n}\Theta\psi_{n})]\\
&=\mathcal{S}_{A} [(f_{n}\Theta g_{n})]/\mathcal{C}_{A} [(\varphi_{n}\Theta\psi_{n})]\\
&=\frac{\mathcal{S}_{A}(f_{n})\mathcal{S}_{A}(g_{n})}{\mathcal{C}_{A}(\varphi_{n})\mathcal{C}_{A}(\psi_{n})}=\frac{\mathcal{S}_{A}(f_{n})}{\mathcal{C}_{A}(\varphi_{n})}\frac{\mathcal{S}_{A}(g_{n})}{\mathcal{C}_{A}(\psi_{n})}\\
&=\mathcal{S}_{A}[f_{n}/\varphi_{n}]\mathcal{S}_{A}[g_{n}/\psi_{n}]=\mathcal{S}_{A}(F)\mathcal{S}_{A}(G).
\end{align*}
(3) Scaling property,  for $k\in\mathbb{R}$ we have
\begin{align*}
&\mathcal{S}_{A}[F(kt)](u)=\frac{1}{\sqrt{2\pi ib}}e^{\frac{i}{2}(\frac{d}{b})s^{2}}\int_{-\infty}^{\infty}sin(\frac{s}{b}t)e^{\frac{i}{2}\frac{a}{b}t^{2}}F(kt)dt\\
&=\frac{1}{k\sqrt{2\pi ib}}e^{\frac{i}{2}(\frac{d}{b})s^{2}}\int_{-\infty}^{\infty}sin(\frac{s}{kb}x)e^{\frac{i}{2}\frac{a}{b}(\frac{x}{k})^{2}}F(x)dx\\
&=\frac{e^{(1-\frac{1}{k^{2}})\frac{i}{2}\frac{d}{b}s^{2}}}{k\sqrt{2\pi ib}}e^{\frac{i}{2}(\frac{d}{b})(\frac{s}{k})^{2}}\int_{-\infty}^{\infty}sin\bigg{(}\frac{s}{k}\frac{1}{b}x\bigg{)}e^{\frac{i}{2}\frac{a}{b}x^{2}}e^{(\frac{1}{k^{2}}-1)\frac{i}{2}\frac{a}{b}x^{2}}F(x)dx\\
&=\frac{e^{(1-\frac{1}{k^{2}})\frac{i}{2}\frac{d}{b}s^{2}}}{k\sqrt{2\pi ib}}e^{\frac{i}{2}(\frac{d}{b})(\frac{s}{k})^{2}}\int_{-\infty}^{\infty}sin\bigg{(}\frac{s}{k}\frac{1}{b}t\bigg{)}e^{\frac{i}{2}\frac{a}{b}t^{2}}e^{(\frac{1}{k^{2}}-1)\frac{i}{2}\frac{a}{b}t^{2}}F(t)dt\\
&=\frac{1}{k}e^{(1-\frac{1}{k^{2}})\frac{i}{2}\frac{d}{b}s^{2}}\mathcal{S}_{A} (e^{(\frac{1}{k^{2}}-1)\frac{i}{2}\frac{a}{b}t^{2}}F)(\frac{s}{k}).
\end{align*}
(4) Shifting property,  for $\tau\in\mathbb{R}$ we have
\begin{align*}
&\mathcal{S}_{A} [F (t+\tau)](u)=\frac{1}{\sqrt{2\pi ib}}e^{\frac{i}{2}(\frac{d}{b})s^{2}}\int_{-\infty}^{\infty}sin(\frac{s}{b}t)e^{\frac{i}{2}\frac{a}{b}t^{2}}F(t+\tau)dt\\
&=\frac{1}{\sqrt{2\pi ib}}e^{\frac{i}{2}(\frac{d}{b})s^{2}}\int_{-\infty}^{\infty}sin(\frac{s}{b}(x-\tau))e^{\frac{i}{2}\frac{a}{b}(x-\tau)^{2}}F(x)dx\\
&=\frac{e^{\frac{i}{2}(\frac{d}{b})s^{2}}}{\sqrt{2\pi ib}}\int_{-\infty}^{\infty}e^{\frac{i}{2}\frac{a}{b}(x-\tau)^{2}}[sin(\frac{s}{b}x)cos(\frac{s}{b}\tau)-cos(\frac{s}{b}x)sin(\frac{s}{b}\tau)]F(x)dx\\
&=\frac{e^{\frac{i}{2}(\frac{d}{b})s^{2}}}{\sqrt{2\pi ib}}\int_{-\infty}^{\infty}e^{\frac{i}{2}\frac{a}{b}(t-\tau)^{2}}[sin(\frac{s}{b}t)cos(\frac{s}{b}\tau)-cos(\frac{s}{b}t)sin(\frac{s}{b}\tau)]F(t)dt\\
&=e^{\frac{i}{2}(\frac{a}{b})\tau^{2}}cos(\frac{s}{b}\tau)\mathcal{S}_{A}(e^{-it\tau(\frac{a}{b})} F)(s) - e^{\frac{i}{2}(\frac{a}{b})\tau^{2}} sin(\frac{s}{b}\tau)\mathcal{C}_{A}(e^{-it\tau(\frac{a}{b})} F)(s).
\end{align*}
 (5) \begin{align*}
 &\mathcal{S}_{A}(e^{ixt} F)(s)=\frac{1}{\sqrt{2\pi ib}}e^{\frac{i}{2}(\frac{d}{b})s^{2}}\int_{-\infty}^{\infty}sin(\frac{s}{b}t)e^{\frac{i}{2}\frac{a}{b}t^{2}}e^{ixt} F(t)dt\\
 &=\frac{1}{\sqrt{2\pi ib}}e^{\frac{i}{2}(\frac{d}{b})s^{2}}\int_{-\infty}^{\infty}sin(\frac{s}{b}t)e^{\frac{i}{2}\frac{a}{b}t^{2}}[cos(xt)+isin(xt)] F(t)dt\\
 &=\mathcal{S}_{A}(cos(xt) F)(s)+i \mathcal{S}_{A} (sin(xt) F)(s).
\end{align*}

(6)
\begin{align*}
&\mathcal{S}_{A}(cos(xt) F(t))(s)=\frac{1}{\sqrt{2\pi ib}}e^{\frac{i}{2}(\frac{d}{b})s^{2}}\int_{-\infty}^{\infty}sin(\frac{s}{b}t)e^{\frac{i}{2}\frac{a}{b}t^{2}}cos(xt) F(t)dt\\
&=\frac{1}{2\sqrt{2\pi ib}}e^{\frac{i}{2}(\frac{d}{b})s^{2}}\int_{-\infty}^{\infty}e^{\frac{i}{2}\frac{a}{b}t^{2}}[sin(\frac{(s+bx)}{b}t)+sin(\frac{(s-bx)}{b}t)]F(t)dt\\
&=\frac{e^{\frac{i}{2}(\frac{d}{b})(bx)^{2}+2bx}}{2\sqrt{2\pi ib}}e^{\frac{i}{2}(\frac{d}{b})(s+bx)^{2}}\int_{-\infty}^{\infty}e^{\frac{i}{2}\frac{a}{b}t^{2}}sin(\frac{(s+bx)}{b}t)F(t)dt\\
&\quad+\frac{e^{\frac{i}{2}(\frac{d}{b})(bx)^{2}-2bx}}{2\sqrt{2\pi ib}}e^{\frac{i}{2}(\frac{d}{b})(s-bx)^{2}}\int_{-\infty}^{\infty}e^{\frac{i}{2}\frac{a}{b}t^{2}}sin(\frac{(s-bx)}{b}t)F(t)dt\\
&=\frac{e^{\frac{i}{2}(db)x^{2}}}{2}\big{[}e^{idxs}\mathcal{S}_{A} \big{(}F\big{)}(s+bx)+e^{-idxs}\mathcal{S}_{A} \big{(}F\big{)}(s-bx)\big{]}.
\end{align*}
(7)
\begin{align*}
&\mathcal{S}_{A}(sin(xt) F(t))(s)=\frac{1}{\sqrt{2\pi ib}}e^{\frac{i}{2}(\frac{d}{b})s^{2}}\int_{-\infty}^{\infty}sin(\frac{s}{b}t)e^{\frac{i}{2}\frac{a}{b}t^{2}}sin(xt) F(t)dt\\
&=\frac{1}{2\sqrt{2\pi ib}}e^{\frac{i}{2}(\frac{d}{b})s^{2}}\int_{-\infty}^{\infty}e^{\frac{i}{2}\frac{a}{b}t^{2}}[cos(\frac{(s-bx)}{b}t)-cos(\frac{(s+bx)}{b}t)]F(t)dt\\
&=\frac{e^{\frac{i}{2}(\frac{d}{b})(bx)^{2}-2bx}}{2\sqrt{2\pi ib}}e^{\frac{i}{2}(\frac{d}{b})(s-bx)^{2}}\int_{-\infty}^{\infty}e^{\frac{i}{2}\frac{a}{b}t^{2}}cos(\frac{(s+bx)}{b}t)F(t)dt\\
&\quad-\frac{e^{\frac{i}{2}(\frac{d}{b})(bx)^{2}+2bx}}{2\sqrt{2\pi ib}}e^{\frac{i}{2}(\frac{d}{b})(s+bx)^{2}}\int_{-\infty}^{\infty}e^{\frac{i}{2}\frac{a}{b}t^{2}}cos(\frac{(s-bx)}{b}t)F(t)dt\\
&=\frac{e^{\frac{i}{2}(db)x^{2}}}{2}\big{[}e^{-idxs}\mathcal{C}_{A} \big{(}F\big{)}(s-bx)-e^{idxs}\mathcal{C}_{A} \big{(}F\big{)}(s+bx)\big{]}.
\end{align*}$\hfill\blacksquare$
\begin{theorem}
If $\delta-\lim F_{n}=F$, then $\mathcal{S}_{A} (F_{n})\rightarrow \mathcal{S}_{A} (F)$ as $n\longrightarrow \infty$ on every compact set of $\mathbb{R}$.
\end{theorem}
\textbf{Proof:} Let $(\delta_{m})$ be a delta sequence such that $F_{n}\Theta\delta_{m}, F\Theta\delta_{m}\in \mathcal{L}^{1}(\mathbb{R}) $ for all $n,m\in \mathbb{N}$ and $\|((F_{n}-F)\Theta\delta_{m})\|_{2}\rightarrow 0$ as $n\rightarrow \infty$ for each $m\in\mathbb{N}$. Let $M$ be a compact subset of $\mathbb{R}$ then $\mathcal{C}_{A}(\delta_{m})> 0$ on $M$ for almost all $m\in\mathbb{N}$. Since $\mathcal{C}_{A}(\delta_{m})$ is a continuous function and $\mathcal{S}_{A}(F_{n})\Theta \mathcal{S}_{A}(\delta_{m})- \mathcal{S}_{A}(F)\Theta \mathcal{C}_{A}(\delta_{m})=\big{(}(\mathcal{S}_{A}(F_{n})- \mathcal{S}_{A}(F))\big{)}\Theta \mathcal{C}_{A}(\delta_{m}))$, implies $\|\big{(}(\mathcal{S}_{A}(F_{n})- \mathcal{S}_{A}(F))\Theta \mathcal{C}_{A}(\delta_{m})\big{)}\|_{1}\rightarrow 0$ as $n\rightarrow \infty$, thus $\mathcal{S}_{A}(F_{n})\rightarrow \mathcal{S}_{A}(F)$ as $n\longrightarrow \infty$ on each $M$.$\hfill\blacksquare$
 

\end{document}